\documentclass[a4paper, 12pt, reqno]{amsart} 
\usepackage{transparent} 
\usepackage{enumitem}
\usepackage[colorinlistoftodos]{todonotes} 
\usepackage{amsmath, amsthm,bbm, amssymb, amsfonts, fullpage,hyperref} 
\usepackage[maxbibnames=99,backend=biber,bibencoding=utf8,style=numeric,giveninits=true,url=false,doi=false,isbn=false]{biblatex} 
\usepackage{mathtools} 
\usepackage{mathrsfs}
\usepackage{amsmath}
\mathtoolsset{showonlyrefs} 
\usepackage{graphicx, subfigure, color,calc} 
\numberwithin{equation}{section} 
\date{\today}

\newcommand{\R}{\mathbb{R}} 
 
\newcommand{\N}{\mathbb N} 
\newcommand{\E}{\mathbb{E}} 
 
\theoremstyle{plain} 
\newtheorem{thm}{Theorem}[section]

\newtheorem{lemma}[thm]{Lemma}

\newtheorem{prop}[thm]{Proposition} 
\theoremstyle{remark} 
\newtheorem{remark}[thm]{Remark} 
\newtheorem{examples}[thm]{Examples}

\newcommand{\lb}{\left(} 
\newcommand{\rb}{\right)} 
\renewcommand{\l}{\lambda} 
\renewcommand{\O}{\mathcal O} 
\newcommand{\lv}{\lvert} 
\newcommand{\rv}{\rvert} 
\newcommand{\lvv}{\left\lvert} 
\newcommand{\rvv}{\right\rvert} 
 
\renewcommand{\t}{\tau} 
\newcommand{\g}{\gamma} 
\newcommand{\s}{\sigma} 
 
\renewcommand{\b}{\beta} 
\renewcommand{\epsilon}{\varepsilon} 
\newcommand{\e}{\epsilon}

\renewcommand{\hat}{\widehat} 
\renewcommand{\P}{\mathbb P} 
\renewcommand{\d}{\delta} 
\newcommand{\supp}{\operatorname{supp}} 
 
\newcommand{\mc}{\mathcal} 
\newcommand{\n}{{(n)}}

\newcommand{\Var}{\operatorname{Var}} 
\newcommand{\Cov}{\operatorname{Cov}} 
\newcommand{\dopp} [1] {\mathbbmss{#1}} 
\newcommand{\SC}{\textup{SC}}

\title{A Central Limit Theorem for Random Variables under a Constraint} 

\author{Jack Cromwell} 
\author{Martin Venker} 
\email{jack.cromwell3@mail.dcu.ie, martin.venker@dcu.ie}
\address{School of Mathematical Sciences, Dublin City University, Dublin, Ireland}

\begin{document}

\begin{abstract}
	We consider vectors chosen randomly from high-dimensional hypersurfaces of the form $\{x\in\R^n: \sum_{j=1}^n V_j(x_j)=an\}$. These in general asymmetric surfaces generalize classical spheres and include unbounded and disconnected surfaces. We show that under mild assumptions as $n\to\infty$, any weighted sum of the random coordinates is asymptotically Gaussian, provided that no weight dominates. An explicit reduction in the variance compared to the unconstrained case is observed. We also study a model with penalization of deviations from the constraint, interpolating between the unconstrained and the constrained case. Our approach involves an explicit representation of the characteristic function of the constrained random variables in terms of independent random variables. 
\end{abstract}

	\maketitle 
	
	\section{Introduction and Statement of Main Results} 
	The central limit theorem (CLT) is generally regarded as one of the most important results in probability theory. In its simplest form it states that for a sequence $X_1,X_2,\dots$ of independent and identically distributed (i.i.d.) real-valued random variables with common mean $\mu\in\R$ and finite variance $\s^2>0$ we have 
	\begin{align} 
		\frac{1}{\sqrt {\s^2n}}\sum_{j=1}^n(X_j-\mu)\to Z\label{CLT} 
	\end{align} 
	in distribution as $n\to\infty$, where $Z$ is a random variable with a standard normal distribution. While it is standard to extend this classical CLT to independent but non-identically distributed random variables, extending to dependent random variables generally depends on the type of dependence.

	In this paper we study sequences of random variables under explicit dependencies introduced by the geometrical constraint 
	\begin{align} 
		\sum_{j=1}^n V_j(X_j)=an\label{constraint} 
	\end{align} 
	for some $a>0$ and functions $V_j$, as $n\to\infty$. A classical example is the $l_2$-sphere in $\R^n$, obtained for functions $V_j(x_j):=x_j^2$ for all $j$, see \cite{DiaFr} for a historical overview. More generally, $V_j(x_j) := |x_j|^p$ for $0<p<\infty$ gives the $l_p$-spheres. Much work has been done on $l_p$-spheres, balls and related spaces as we will outline below. This is usually achieved by explicit representations of distributions on these balls or spheres obtained by normalizing a vector of i.i.d random variables. For example, the uniform distribution on the $l_2$ unit sphere can be realized by normalizing a vector of i.i.d.~Gaussians, and the cone measure on the $l_p$-sphere arises in the same way from i.i.d generalized Gaussian variables. Such representations rely on the homogeneity of $x_j \to |x_j|^p$ as well as the identical distribution of all coordinates, and are not available for general $V_j$. A main motivation of this paper is to develop an approach that does not require homogeneity, symmetry, or identically distributed coordinates.

    Allowing the $V_j$ to differ with $j$, and even change with $n$, leads to a wide variety of hypersurfaces. For example, $V_j(x) := c_jx^2$ with different $c_j >0$ gives ellipsoids. Bounded $V_j$ can lead to unbounded surfaces and $S_{a,V}$ may consist of disconnected pieces, with the number of pieces growing with $n$. Examples will be discussed in more detail below, see Examples \ref{Examples}. Our assumptions on the $V_j$ will be very general and mild, generally relying on a small bit of curvature of $S_{a,V}$. The standard simplex, obtained for $V_j(x) : = x$ on $[0,\infty)$, lies at the boundary of our framework and is excluded by our curvature requirements, however, for this case we point the reader to \cite{GamboaVenker} where random vectors under several linear constraints are studied.

	To make our setup precise, let for each $n\geq1$ continuous functions $V_1^\n,\dots,V_n^\n:\R\to[0,\infty)$ and $a^\n>0$ be given and define 
	\begin{align} 
		S_{a,V}:=\left\{x\in\R^n:\sum_{j=1}^nV_j^\n(x_j)=a^\n n\right\}.\label{S} 
	\end{align} 
	
	The superscript $^\n$ signifies that we can  change $a^\n$ and the functions $V_j^\n$ completely going from $n$ to $n+1$, but we will usually suppress this dependence in the notation for improved readability, simply writing $a$ and $V_j$. 
	As $S_{a,V}$ is a hypersurface in $\R^n$, it is not obvious which measure to equip it with. Even when $S_{a,V}$ is compact, there are several natural candidates for a uniform distribution on it. This relates to the classical fact that conditioning a probability measure on a null set is not well-defined unless the way of conditioning is specified (Borel-Kolmogorov paradox). We will equip $S_{a,V}$ with the so-called  \textit{microcanonical measure}, which is the most natural measure from a statistical mechanics point of view. To define it, let $H^{n-1}$ denote the $n-1$ dimensional Hausdorff measure in $\R^n$. We will assume that $H^{n-1}(S_{a,V})>0$ going forward. Setting 
	\begin{align}
		V:\R^n\to\R,\quad V(x):=\sum_{j=1}^nV_j(x_j),\label{def:V}
	\end{align}
	the microcanonical measure is defined as
	\begin{align}
		d\mu_{S_{a,V}}(x):=\frac{1}{\|\nabla V(x)\|}dH^{n-1}(x),\label{def:mu}
	\end{align} where, under the regularity assumptions below $\nabla V(x)$ exists for $H^{n-1}$ almost every $x$, see Remark \ref{R2}(4).  The name stems from its occurrence in microcanonical ensembles of statistical physics. The following remark addresses the choice of measure on $S_{a,V}$.
\begin{remark}[On different measures on $S_{a,V}$] \label{remark_measures} There are at least three different measures that have been considered in the literature as measures for surfaces like $S_{a,V}$, the Hausdorff or surface measure $H^{n-1}$, the microcanonical measure $\mu_{S_{a,V}}$, and the cone measure.

    While the Hausdorff measure might seem the most natural from a geometric point of view, 
	the co-area formula from geometric measure theory (\cite[Theorem 3.2.12]{Federer}) shows
	\begin{align}
		\int_{\R^n}g(x)dx=n\int_{\R}\lb\int_{S_{a,V}}g(x)\frac{1}{\|\nabla V(x)\|}dH^{n-1}(x)\rb da=n\int_{\R}\lb\int_{S_{a,V}}g(x)d\mu_{S_{a,V}}(x)\rb da\label{coarea}
	\end{align}
	with $\mu_{S_{a,V}}$ of \eqref{def:mu}, and valid for any $g\in L^1(\R^n)$ under mild assumptions on the $V_j$. Formula \eqref{coarea} establishes the microcanonical measure $\mu_{S_{a,V}}$ as the natural measure when disintegrating via level sets. In particular, if $X$ has a density $g$, then the conditional law of $X$ given $V(X) = an$ defined by this disintegration is proportional to $g d \mu_{S_{a,V}}$. Another instructive view on the difference of $H^{n-1}$ and $\mu_{S_{a,V}}$ is looking at how they derive from the Lebesgue measure in ambient space. The microcanonical measure is obtained from scaling the $n$-dimensional Lebesgue measure of the thin shell $\{x\in\R^n: an-\e\leq\sum_{j=1}^n V_j(x_j)\leq an+\e\}$ as $\e\to0$. We will see a variant of this in Proposition \ref{prop:convergence} below where $\mu_{S_{a,V}}$ arises when approaching the constraint $\sum_{j=1}^nV_j(X_j)=an$ using the canonical Gaussian approximation of the ``delta function" $\d(\sum_{j=1}^nV_j(x_j)-an)$. In contrast to that, the Hausdorff measure $H^{n-1}$ can be obtained using a non-Gaussian approximation of the delta function, obtained by weighting the Gaussian approximation of the delta function by the factor $\|\nabla V(x)\|$. Thus from the point of conditioning the Lebesgue measure to $S_{a,V}$, the microcanonical measure is more natural. It should be noted that in the classical examples of the $l_2$-sphere and the simplex (part of the $l_1$-sphere), Hausdorff and microcanonical measure are multiples of each other as $\|\nabla V(x)\|$ is constant on these surfaces. 
	
	Finally, the cone measure of a measurable subset $A$ of $S_{a,V}$ is defined as the Lebesgue measure of the cone $\{tx:x\in A,\, t\in[0,1]\}$. To define this measure, $S_{a,V}$ must be star-shaped, which is in general not the case. Where it is defined, it has a density proportional to 
    \begin{align}
        \frac{\langle x,\nabla V(x)\rangle}{\|\nabla V(x)\|}dH^{n-1}(x).
    \end{align}
    For $l_p$-spheres, where $V_j(x_j)=\lv x_j\rv^p$, Euler's identity for homogeneous functions gives $\langle x, \nabla V(x) \rangle = p V(x) = pan$ on $S_{a,V}$, so the cone measure and microcanonical measure are multiples of each other. For more general functions $V_j$ however, these two measures differ. Let us finally remark that the on the $l_p$-sphere the difference between the cone measure, and thus the microcanonical measure, and the Hausdorff measure, was shown in \cite{NaorRomik,Naor} to be small as $n\to\infty$.
\end{remark}
	
	The microcanonical measure $\mu_{S_{a,V}}$ may be finite or infinite. If it is finite, the normalized version of it may be considered as uniform distribution on $S_{a,V}$. We note in passing that we  do not exclude points of $S_{a,V}$ at which $\nabla V$ vanishes. At such points the density of $\mu_{S_{a,V}}$ with respect to $H^{n-1}$ is singular, but our assumptions below will guarantee that the measures considered below are finite for sufficiently large $n$.
    
	 Let now measurable functions $f_j=f_j^\n:\R\to[0,\infty)$ be given and define 
	\begin{align} 
		dP_{n}(x):=\frac{1}{Z_{n}}\prod_{j=1}^nf_j(x_j)d\mu_{S_{a,V}}(x), \quad Z_{n}:=\int_{S_{a,V}}\prod_{j=1}^nf_j(x_j)d\mu_{S_{a,V}}(x) \label{P_n}.
	\end{align} 
	$P_n$ is a probability measure on $S_{a,V}$, where $Z_n$ is its normalizing constant. We assume here that $Z_{n}>0$, finiteness of $Z_n$ for large $n$ will follow from Proposition \ref{prop:convergence} below. The reader might want to think of the probability measure $P_{n}$ as realizing random variables $X_j$ with densities proportional to $f_j$ but subject to the constraint \eqref{constraint}. While the functions $V_j$ will be continuous, we will allow for the $f_j$ to have a finite number of discontinuities, hence allowing us to restrict to different regions of $S_{a,V}$ using indicator functions. We observe that under the probability measure $P_n$, the coordinate random variables $X_1,\dots,X_n$ are not independent and neither are they identically distributed.
	
	The constraint \eqref{constraint} does not only pose the technical challenge of rendering the support of the probability measure lower-dimensional, it also in general strongly influences the behavior of the random variables $X_1,\dots,X_n$. For example, it may force otherwise unbounded random variables to be bounded. \\
    
		 While \eqref{constraint} might be understood as a \textit{hard constraint}, we will also consider a model representing a \textit{soft constraint}, under which deviations of $\sum_{j=1}^n V_j(X_j)$ from $an$ are possible but penalized. To this end we introduce, for a parameter $\g>0$ the probability density on $\R^n$ defined as 
	\begin{align} 
		P_{n,\SC}(x):=\frac{1}{Z_{n,\SC}}e^{-\g\lb \frac1{\sqrt n}\sum_{j=1}^n V_j(x_j)-a\sqrt n\rb^2}\prod_{j=1}^nf_j(x_j)\label{P_n,SC} 
	\end{align} 
	with normalizing constant $Z_{n,\SC}$, the subscript SC standing for ``Soft Constraint''. Finiteness of $Z_{n,\SC}$ is not automatic as integration is now over all of $\R^n$ rather than over $S_{a,V}$. However, under our assumptions below it will hold for every $\g >0$, see \eqref{manipulation} below. We will slightly abuse notation by referring to probability densities like $P_{n,\SC}$ as probability measures when there is no danger of ambiguity. The exponential factor in \eqref{P_n,SC} penalizes deviations from the constraint \eqref{constraint}, the strength of this penalization depending on the value of $\g$. It formally interpolates between the unconstrained case $\g=0$ and the hard constraint obtained for $\g\to\infty$ (Proposition \ref{prop:convergence}). As the analysis of $P_{n,\SC}$ is similar to the one of $P_{n}$, and in fact the analysis of the hard constraint case $P_n$ proceeds via the soft constraint case $P_{n,\SC}$, we will treat both cases in parallel.\\ 
	
	We will compare the limiting behavior of the dependent $X_1,\dots,X_n$ under $P_{n}$ or $P_{n,\SC}$ to the limiting behavior of independent $X_1,\dots,X_n$ under a different probability measure with density 
	\begin{align} 
		P_{n,\b}(x):=\frac{1}{Z_{n,\b}}e^{-\b\sum_{j=1}^n V_j(x_j)}\prod_{j=1}^nf_j(x_j),\label{P_n,b} 
	\end{align} 
	where $\b\geq0$ such that the normalizing constant $Z_{n,\b}$ is finite. In the spirit of classical statistical mechanics, $P_n$ corresponds to a microcanonical ensemble of fixed ``energy" $\sum_{j=1}^nV_j(x_j)$, whereas
	$P_{n,\b}$ corresponds to a canonical ensemble at inverse temperature $\b$.  $P_{n,\SC}$ corresponds to an interpolating ensemble, it has in fact been considered as such in the statistical mechanics literature before,  see \cite{Challa}, and \cite{ACV} for an application in random matrix theory.

    Expectations of local quantities are classically expected to agree between the two ensembles as $n \to \infty$ by equivalence of ensembles. However, this equivalence does not extend to fluctuations, indeed our results rigorously quantify this difference. Our approach will link the two ensembles directly for every finite $n$, not just in the limit. Proposition \ref{prop_CF} below expresses the characteristic functions of $P_n$ and $P_{n,\SC}$ exactly in terms of $P_{n,\b}$, with $\b$ chosen such that the constraint holds in expectation under $P_{n,\b}$. \\

 In order to state our results we now move to assumptions. Our first set of assumptions is about regularity of the functions $f_j$ and $V_j$ but also concerns the curvature of $S_{a,V}$. We remark in advance that the assumptions are fairly general. A main theme is that, as $a$, $f_j$ and $V_j$ may depend on $n$, some form of uniformity in $j$ and $n$ has to be assumed. We will for the sake of brevity, and with a slight abuse of notation, write ``independent of $j$" to imply independence from both $j$ and $n$. 
 \\

\noindent
\textbf{Assumptions 1} (Analytic Assumptions)
\makeatletter
\def\@currentlabel{1}\label{ASS}
\makeatother
\begin{enumerate}[label=\textbf{1.\arabic*}, ref=1.\arabic*]
    \item \textbf{(Regularity)} \label{ass:regularity}
    There are $M>0$ and $p,T \in \{2,3,\dots\}$, independent of $j$, such that for each $j$ the following holds 
    \begin{itemize}
        \item $f_j: \mathbb{R} \to \mathbb{R}_+$ is continuously differentiable except for a finite set of points $D_j$ at which jump discontinuities may occur. Furthermore, $f_j(x) \leq M$, for all $x\in\R$.
        \item $V_j : \R \to \R_+$ is continuous on $\R$, and $p$ times continuously differentiable on $\R \setminus H_j$ for a finite set $H_j$.
        \item $|D_j| + |H_j| \leq T$.
    \end{itemize}
    \item \textbf{(Point of non-vanishing curvature)} \label{ass:curve}
    There exist constants $M>0$, $\d>0$ and $\e>0$, independent of $j$ such that for each $j$ there is an interval $J_j \subset \R \setminus H_j$ of length $|J_j| > \d$ and 
    \begin{align}
        |V''_j(x)| > \e, \quad V_j(x) \leq M, \quad f_j(x) > \e, \quad \text{ for all } x \in J_j.
    \end{align}
    \item \textbf{(Non-flatness for sufficiently many $j$)} \label{ass:p}
    There exist $c>0$, $T \in \N$ and $\e>0$ independent of $j$ such that for each $n$, there is a set $G =G^\n \subset \{1,\dots,n\}$ with $|G|>cn$ for which the following holds for all $j \in G$: There exist $A_j \leq T$ intervals $I_{j,1},\dots,I_{j,A_j}$ with
    \begin{align}
        \lb \supp (f_j) \setminus H_j\rb \subset \bigcup_{k=1}^{A_j} I_{j,k},
    \end{align}
    and $p_{j,k} \in \{1,\dots,p\}$, with $p$ as in Assumption \ref{ass:regularity},  such that
     \begin{align}
         |V_j^{(p_{j,k})}(x)| \geq \e \quad  \text{ for all } x \in I_{j,k} \text{ and for all } k \in \{1,\dots,A_{j}\},
     \end{align}
     with $p_{j,k}=1$ whenever $I_{j,k}$ is unbounded.
    \item \textbf{(Subcritical mass)} \label{ass:mass}
    There are constants $\e>0$ and $\d>0$ independent of $j$ such that for each $j$ there exists an interval $I_j$ with length $|I_j| >\d$ and
    \begin{align}
        f_j(x) \geq \e \text{ and } V_j(x) \leq a-\e \text{ on } I_j.
    \end{align}
   
\end{enumerate}
\begin{remark} \label{R2}
	\begin{enumerate}
		\item As mentioned above, the generality of the setup with all relevant quantities depending on $j$ and $n$, requires uniformity in the assumptions. For the simplified, yet general setup of $V_j=V_j^\n=\bar V,\,f_j=f_j^\n=\bar f,\,a^\n=a$ for all $j$ and $n$ for fixed functions $\bar V$ and $\bar f$, and fixed $a>0$, all uniformity assumptions become trivial. Assumption \ref{ass:curve} holds as soon as $\bar V$ is not affine on some interval on which $\bar f$ is bounded away from 0. Assumption \ref{ass:p} holds with $G=\{1,\dots,n\}$ if at every point of $\supp(\bar f)$ some derivative $\bar V^{(k)}$ is nonzero, and if $\supp(\bar f)$ is unbounded, $|\bar V'|$ is bounded away from zero on $\supp(\bar f)$ outside of some bounded interval. 
        Finally, Assumption \ref{ass:mass} holds as soon as $\bar V<a$ at some point of an open interval on which $\bar f$ is bounded away from 0.
	\item As will become clear below, a major element of our approach is a representation of the characteristic function of $P_n$ in terms of independent random variables, see Proposition \ref{prop_CF}. The analysis requires uniform decay of the one-dimensional Fourier transforms
    \begin{align}
        \int_\R e^{irV_j(x)}f_j(x)e^{-\b V_j(x)}dx \text{ as } r \to \infty.
    \end{align}
    Assumption \ref{ass:p} is designed to provide that decay at the rate $|r|^{-1/p}$ by means of the van der Corput lemma, see Lemma \ref{lem_vdc}. Without this assumption, $V_j$ could be constant on some interval of positive length, allowing for $S_{a,V}$ to have pieces which are flat in the $x_j$ direction. Thus $P_n$ may put mass on slices $\{x_i = z\}$, meaning coordinates $X_j$ may have atoms. This assumption prevents such atoms. It however allows $V_j$ to have degenerate stationary points of bounded order, which only affect the rate through $p$.  The assumption is only necessary for a positive proportion of the coordinates so the remaining $V_j$ are only subject to Assumptions \ref{ass:regularity}, \ref{ass:curve}, \ref{ass:mass}. Assumption \ref{ass:p} can be relaxed to allow for weaker or changing types of decay; this would however lengthen the technical parts of this paper significantly and obstruct the main new ideas.
    \item Assumption \ref{ass:mass} is only used to guarantee the existence of a $\b$ under which the constraint holds in expectation.
    \item Assumption \ref{ass:regularity} allows each $V_j$ to be non-differentiable at finitely many points $H_j$, at which $V_j'$ may blow up or have a jump. Examples are $|x|^r$ for $0<r<1$. Consequently, $\nabla V(x)$ exists and is finite only for $x$ outside the set $\{x: x_j \in H_j \text{ for some } j\}$. On $S_{a,V} \cap \supp (\prod_{j=1}^n f_j(x_j))$, this exceptional set has $H^{n-1}$ measure zero (see proof of Proposition \ref{prop:convergence}), and hence does not affect $P_n$. We note that in the case of $|x|^r$ for $r < 1$, the integrability conditions in Assumption \ref{ass:int} become significantly stronger and we can no longer utilize $f_j(x) :=1$.
	\end{enumerate}
	
\end{remark}
\medskip 
With our analytic setup in place, we move to two more assumptions of a more probabilistic nature.
Our main results provided in Theorem $\ref{T1}$ and Theorem $\ref{T2}$ below involve applying a CLT to a sequence of non-identically distributed random variables given by linear combinations of $V_j(X_j)$ and $\sqrt n \l_j X_j$, where $X:=(X_1,\dots,X_n)$ has the distribution $P_{n,\b}$. As with a classical CLT we make moment assumptions to ensure validity of our results. One problem with this is however that the $\b$, while unique, is implicit as it is determined by the $V_j,\,f_j$ and $a$, and hence assumptions using $P_{n,\b}$ would be hard to check. We overcome this by making assumptions w.r.t.~$P_{n,\b_0}$ for some fixed $\b_0\geq0$. We will show later on that the corresponding assertions for $P_{n,\b}$ will follow from the ones for $P_{n,\b_0}$. Let $\E_{{n,\b}}$, $\Var_{{n,\b}}$ and $\Cov_{{n,\b}}$ denote expectation, variance and covariance with respect to $P_{n,\b}$ of $\eqref{P_n,b}$. 

\vspace{1em}
\noindent
\noindent
\textbf{Assumptions 2} (Moment Assumptions)
\makeatletter
\def\@currentlabel{2}\label{MASS}
\makeatother
There exist constants $\b_0\in[0,\infty)$ and $K>0$ such that the following holds.
\begin{enumerate}[label=\textbf{2.\arabic*}, ref=2.\arabic*]
	 \item \textbf{(Integrability)} \label{ass:int}We have for each $j$
	\begin{align}
		&\int_\R f_j(x) e^{-\b_0 V_j(x)} dx\le K \label{Zbound}\quad\text{ and } \\
		&\int_{\R\setminus (D_j \cup H_j)} \lvv f'_j(x)e^{-\b_0 V_j(x)} \rvv + \lvv f_j(x) V'_j(x) e^{-\b_0V_j(x)}\rvv +\lvv f_j(x) V''_j(x) e^{-\b_0V_j(x)}\rvv dx \leq K.\label{BI}
	\end{align}
    \item \textbf{(Moment Bounds)} \label{ass:BM} We have for each $j$ and each $n$
		\begin{align} 
				&\E_{{n,\b_0}}(V_j(X_j)^3) \leq K,\quad \E_{{n,\b_0}}(|X_j|^3) \leq K \text{ and}\label{moment_existence}\\
				&\E_{{n,\b_0}}(\sum_{j=1}^n V_j(X_j))\geq na. \label{Mbound}
		\end{align} 
         For the soft constraint ensemble $P_{n,\SC}$ with parameter $\g> 0$, we will replace condition \eqref{Mbound} with the slightly stronger assumption
         \begin{align}
         	\E_{{n,\b_0}}(\sum_{j=1}^n V_j(X_j))\geq n(a + \frac{\b_0}{2\g}) \label{Mbound2}.
         \end{align} 
\end{enumerate}

\begin{remark}
The constant $\b_0$ in Assumptions 2 acts as a lower bound of the actual $\b$ we will work with. In the hard constraint case, $\b$ is chosen such that the constraint holds in expectation, i.e $\E_{n,\b} \lb \sum_{j=1}^n V_j(X_j) \rb = an$, and Assumption \ref{ass:mass} paired with condition \eqref{Mbound} guarantees such a $\b \geq \b_0$ exists. Similarly, in the soft constraint case, \eqref{Mbound2} ensures the existence of a $\b \geq \b_0$ satisfying the modified constraint $\E_{n,\b} \lb \sum_{j=1}^n V_j(X_j) \rb = \lb a+ \b /(2\g) \rb n$. See Proposition \ref{lemma_beta} for further details. 
\end{remark}

    We now state our main results. We will study a generalization of \eqref{CLT}: given a weight vector $\l^\n=\l=(\l_1,\dots,\l_n)$, we consider the statistic 
	\begin{align} 
		\sum_{j=1}^n\l_j(X_j-\E_{{n,\b}}(X_j))=\langle\l,X-\E_{{n,\b}}(X)\rangle,\label{statistic} 
	\end{align} 
	where $\E_{{n,\b}}$ denotes expectation w.r.t.~$P_{n,\b}$ and $\E_{{n,\b}}(X)$  the mean vector of $X$ w.r.t.~$P_{n,\b}$, assuming existence of the latter. While the l.h.s.~of \eqref{statistic} is a typical probabilistic quantity, the r.h.s.~of \eqref{statistic} represents a more geometric point of view. Indeed, if $\|\l\|_2=1$, then \eqref{statistic} is the scalar projection of the centered random vector onto $\l$. From this point of view we study the high-dimensional $S_{a,V}$ by looking at one-dimensional projections.  Our main results below will in particular show that these projections behave similarly as $n\to\infty$ for every direction $\l$ in which no single weight dominates.\\

     Our first result deals with the hard constraint. To state it, we recall that $\|\cdot\|_2$ denotes the Euclidean norm, and we denote by $\|\cdot\|_{\max}$ the maximum norm. Moreover, let $ \mc N(0,1)$ denote the standard normal distribution.
	\begin{thm} $\label{T1}$
    Let  Assumptions 1 and 2 hold and let for each $n$, $\l=\l^\n \in \R^n$ be given such that
    \begin{align}
		\lim_{n\to\infty} \|\l\|_{\max} =0 ,\quad 0< \liminf_{n\to\infty}\|\l\|_2 \leq \limsup_{n\to\infty}\|\l\|_2  < \infty. \label{LASS}
		\end{align}
   There exists a unique $\b\geq \b_0$, bounded in $n$, such that $\E_{n,\b}(\sum_{j=1}^n V_j(X_j)) = na$. With this unique $\b$, we define 
		\begin{align}
			\sigma_n^2 := {\sum_{j=1}^{n}{\l_j}^2{\Var_{{n,\b}}(X_j) -\frac{\left(\sum_{j=1}^{n}\l_j\Cov_{{n,\b}}(X_j,V_j(X_j))\right)^2}{\sum_{j=1}^{n} \Var_{{n,\b}}(V_j(X_j)) }}}.\label{variance}
		\end{align}
      Let $X$ have the distribution $P_{n}$, then as $n \to \infty$ 
        \begin{align}
            \frac{1}{\sigma_n}\sum_{j=1}^n \l_j \lb X_j -\E_{n,\b}(X_j) \rb \to \mc N(0,1)\quad\text{ in distribution.} 
        \end{align}
	\end{thm}
	The next theorem deals with the soft constraint case.
	\begin{thm} $\label{T2}$
    Let  Assumptions 1 and 2 hold, let $\g >0$ be fixed and let for each $n$, $\l=\l^\n \in \R^n$ be given such that \eqref{LASS} holds.
   There exists a unique $\b=\b(\g) \geq \b_0$, bounded in $n$, such that $\E_{n,\b}\left(\sum_{j=1}^n V_j(X_j)\right) = n\lb a+ \frac{\b}{2\g}\rb$. With this unique $\b$, we define 
		\begin{align}
			\s_n^2 := {\sum_{j=1}^{n}{\l_j}^2{\Var_{{n,\b}}(X_j) -\frac{\left(\sum_{j=1}^{n}\l_j\Cov_{{n,\b}}(X_j,V_j(X_j))\right)^2}{\sum_{j=1}^{n}\lb \Var_{{n,\b}}(V_j(X_j)) +\frac{1}{2\g} \rb}}}.\label{variance2}
		\end{align}
        Let $X$ have the distribution $P_{n,\SC}$, then as $n \to \infty$ 
        \begin{align}
            \frac{1}{\sigma_n}\sum_{j=1}^n \l_j \lb X_j -\E_{n,\b}(X_j) \rb \to \mc N(0,1)\quad\text{ in distribution.} 
        \end{align}
	\end{thm}

    \begin{remark}\noindent
    	\begin{enumerate}
    		\item An interesting feature in the CLTs Theorems \ref{T1} and \ref{T2} is the forms of the variances $\s_n^2$. The first term $\sum_{j=1}^{n}{\l_j}^2{\Var_{{n,\b}}(X_j)}$ is the variance of $\sum_{j=1}^n\l_jX_j$ for $(X_1,\dots,X_n)$ having distribution $P_{n,\b}$, under which the $X_j$ are independent. We see that under $P_n$ the variance is reduced by a regression term that captures how much of the linear statistic $\sum_{j=1}^n \l_j X_j$ is explained by the constraint $\sum_{j=1}^n V_j(X_j)=an$. The constrained variance $\s_n^2$ in the hard constraint case is the Schur complement of $\sum_{j=1}^n\Var_{n,\b}(V_j(X_j))$ in the covariance matrix of $(\sum_{j=1}^nV_j(X_j),\sum_{j=1}^n\l_jX_j)$, see the proof of Lemma \ref{sigmalemma} below. Thus it is the residual variance of the linear statistic after subtracting its best linear unbiased predictor based on the constraint. Note that thanks to our curvature Assumption \ref{ass:curve} $\sum_{j=1}^n V_j(X_j)$ is not fully linear, and this will essentially lead to $\s_n$ being bounded away from 0. We will prove this below in Lemma \ref{sigmalemma}.
    		
    		An equivalent formula for the variance \eqref{variance} has been found in \cite{LPV} in the context of statistical mechanics, and then later been shown to hold for a larger class of ensembles and functions in \cite{CancriniOlla}. An extension of the variance formula to several linear constraints has been provided in \cite{GamboaVenker}, where it was found to be a certain orthogonal projection. Finally, we note that \eqref{variance2} provides an interpolation between the variance arising for independent random variables under $P_{n,\b}$ for $\g=0$, and the constrained random variables under $P_n$ for $\g\to\infty$.
    		\item The assumptions on $\l$ in \eqref{LASS} guarantee that no particular $\l_j$ dominates the others. This is in general optimal.
            \item Theorem \ref{T1} shows that any projection of the centered random element from the high-dimensional $S_{a,V}$ onto the real line in directions $\l$ with \eqref{LASS} is asymptotically Gaussian.
    	\end{enumerate}     
    \end{remark}

    We will now illustrate the scope of Theorems \ref{T1} and \ref{T2} by several examples, which have been chosen to exhibit different features the model allows, namely generalizing $l_p$-spheres, disconnected supports, unbounded surfaces, $n$-dependence and asymmetry. Where possible we will compute $\b$ and $\s_n$ explicitly. Recall that $\b$ is determined by the requirement that the constraint holds in expectation, that is 
    \begin{align}
        \E_{n,\b} \lb \sum_{j=1}^n V_j(X_j) \rb = na.
    \end{align}
    Such computations are feasible as under $P_{n,\b}$, the coordinates are independent with explicit one dimensional densities proportional to $f_j(x)e^{-\b V_j(x)}$. Hence both the defining equation for $\b$ and the formula for $\s_n^2$ as in \eqref{variance} reduce to computing one dimensional integrals. In cases where these computations are not explicit, they will be feasible numerically.
 
 \medskip
\begin{examples}\label{Examples}
\makeatletter
\def\@currentlabel{1.7}\label{EX}
\makeatother
\noindent
\begin{enumerate}
    
    \item \textbf{(Mixed exponents sphere)}
     Fix $a>0$ and an integer $P \geq 2$. For each $n$, let $p_1,\dots,p_n \in \{2,\dots,P\}$, and set $V_j(x) := |x|^{p_j}$ and $f_j(x):=1$. Then $P_n$ lives on the compact asymmetric surface
    \begin{align}
        \left\{ x \in \R^n : \sum_{j=1}^n \lv x_j\rv^{p_j} = an \right\}.
    \end{align}
    Unless all $p_j = 2$, $||\nabla V||$ is non constant on this surface so $P_n$ is not the normalized surface measure. It is easy to see that the assumptions hold but we note that Assumption \ref{ass:p} holds with $G = \{1,\dots, n\}$ where we take the $p_j$-th derivative on the intervals $[-1,0),(0,1]$ and the first derivative on $(-\infty,-1]$ and $[1,\infty)$, and the second set of assumptions holds with $\b_0 = 1/(Pa)$.  Under $P_{n,\b},$ the $ X_j$ are all independent and have symmetric density proportional to $\exp(-\b |x_j|^{p_j})$ on $\R$, so substituting $ u = |x|^{p_j}$ shows that $|X_j|^{p_j}$ is Gamma distributed with shape $1/p_j$ and rate $\b$. Hence all moments are explicit with
    \begin{align}
        \E_{n,\b} |X_j|^k = \b^{-k/p_j} \Gamma\lb \frac{k+1}{p_j} \rb \lb \Gamma\lb \frac1 {p_j} \rb\rb^{-1}\text{ for all } k \geq 0,
    \end{align}
     with $\Gamma$ denoting the gamma function. This allows us to compute 
    \begin{align}
          \b = \frac{1}{na} \sum_{j=1}^n \frac{1}{p_j}, \quad \s_n^2 = \sum_{j=1}^n \l_j^2 \b^{-2/p_j} \frac{\Gamma(3/p_j)}{\Gamma(1/p_j)},
    \end{align}
    where we note that the covariance term in \eqref{variance} vanishes due to coordinatewise symmetry.
    \newline \indent
    \item \textbf{(Section of mixed exponents sphere)}
    In order to explicit the reduction in variance in the previous example we set $f_j(x) : = \dopp 1_{[0,\infty)}(x)$ restricting $P_n$ to 
    \begin{align}
        \left\{ x \in [0,\infty)^n : \sum_{j=1}^n x_j^{p_j} = na \right\}.
    \end{align}
    Analogous to above, we can use the independence structure of $P_{n,\b}$, relating this to a Gamma distribution, to compute
    \begin{align}
         \b = \frac{1}{na} \sum_{j=1}^n \frac{1}{p_j}
    \end{align}
    and
    \begin{align}
    \s^2_n = \sum_{j=1}^n \l_j^2 \b^{-2/p_j} \lb \frac{\Gamma(3/p_j)}{\Gamma(1/p_j)} - \lb\frac{\Gamma(2/p_j)}{\Gamma(1/p_j)}\rb ^2 \rb - \frac{1}{na\b} \lb \sum_{j=1}^n \l_j \b^{-1/p_j} \frac{\Gamma(2/p_j)}{p_j \Gamma (1/p_j)} \rb^2.
    \end{align}
    We now can see the  reduction in variance via the subtracted term which is nonzero as the density is no longer symmetric.
    \newline \indent
    Formally taking $p_j =1$ for all $j$, giving  the simplex, violates Assumption \ref{ass:curve}, and would also lead to $\s^2_n = 0$ for $\l = n^{-1/2}(1,\dots,1)$. The problem in this case is that both constraint and the statistic $\sum_{j=1}^n \l_j X_j$ are linear and thus the hard constraint annihilates the randomness for particularly chosen $\l_j$. As this case is however well-known, we do not include it in our framework here and refer to \cite{GamboaVenker} instead, where several linear constraints are studied, and the necessary restrictions on the $\l_j$ are detailed.

      \item \textbf{(Double Well Potential)} 
    Let $V_j(x) := (x^2-1)^2$ and $f_j(x) :=1$ for all $j$. Each $V_j$ has two wells at $x = \pm 1$, separated by a barrier at zero of height $V_j(0) =1$. The constraint gives the coordinates a total energy budget of $an$. If $an <1$, the budget is too small for any single coordinate to climb over the barrier, so every coordinate stays in one of its two wells. Hence, $S_{a,V}$ splits into $2^n$ disjoint pieces. Once $an>1$, a coordinate can cross the barrier while the others give up energy and the pieces merge. The assumptions are easily checked with $\b_0 = 1/(4(a+\sqrt a))$. Since $V_j$ and $f_j$ are even and independent of $j$, the covariance term in \eqref{variance} is zero, and the $X_j$ are identically distributed under $P_{n,\b}$. While $\b$ is only given implicitly, by  $\E_{1,\b}((X^2-1)^2) =a$, the variance term is explicit in terms of $\b$. Integration by parts gives $\E_{1,\b}[X_1 V_1'(X_1)] = \b^{-1}$, and since $xV_1'(x) = 4V_1(x) + 4(x^2-1)$, we obtain
    \begin{align}
        \s_n^2 = ||\l||_2^2 \E_{1,\b}(X_j^2) =  \lb 1- a + \frac{1}{4\b}\rb||\l||^2_2.
    \end{align}
    
    We note that a surface which is disconnected for all $n$ can also be easily achieved by $f_j$ that are discontinuous, for example by taking $f_j(x_j) : =  \dopp 1_{ \{ 1 \leq |x_j| \leq 2\}} $, with $V_j$ as above, $a =1$ and $\b_0 = 0$.

    \item \textbf{(Unbounded surface)}
    For odd $j$ let $V_j(x) := x^2$ and $f_j(x) :=1$. For even $j$ let $V_j(x) := x^2/(1+x^2)$ and $f_j(x) := e^{-x^2}$. Since the even-indexed $V_j$ are bounded, $S_{a,V}$ is unbounded. Assumption \ref{ass:p} holds with $G$ being the odd indices, so any $c < 1/2$ works. The remaining coordinates only need Assumptions \ref{ass:regularity}, \ref{ass:curve} and \ref{ass:mass} which hold near zero. Assumption 2 holds for $\b_0 = 1/(4a)$. Since all $V_j$ and $f_j$ are even, the covariance term in \eqref{variance} is zero. For odd $j$, $X_j$ is Gaussian with variance $1/(2\b)$, so we can write
    \begin{align}
        \s_n^2 = \frac{1}{2\b} \sum_{j \text{ odd}}\l_j^2 + \sum_{j \text{ even}} \Var_{n,\b}(X_j) \l_j^2.
    \end{align}
    For even $j$ the distribution of $X_j$ under $P_{n,\b}$ is not a classical distribution so neither $\b$ nor $\Var_{n,\b}(X_j)$ are available in closed form, however these values can easily be approximated numerically.
    \item \textbf{($n$-dependence)}
    Let $g : [0,1] \to [g_{-},g_{+}]$ be continuous with  constants $0 < g_{-} \leq g_{+}$. Set $V_j^\n (x) := g(j/n) x^2$ and $f_j(x) := \dopp 1_{[0,\infty)}(x)$. Then $P_n$ lies on the positive orthant of $S_{a,V}$, an ellipsoid whose axes vary along the coordinates. The assumptions hold with $\b_0 = 1/(2a)$.
    Under $P_{n,\b}$, $X_j$ has density proportional to $\exp \lb-\b g(j/n)x^2\rb$ on $\R_+$, so $X_j$ is distributed as $|Z|/\sqrt {2 \b g(j/n) }$ with $Z \sim \mc N(0,1)$. In particular $\E_{n,\b}(V_j(X_j)) = \frac{1}{2\b}$ for every $j$, so $\b = 1/(2a)$. Using 
    \begin{align}
        E|Z| = \sqrt{2/\pi},\quad E|Z|^3 = 2 \sqrt{2/\pi} \text{ and }\Var(Z^2) =2,
    \end{align}
    we obtain
\begin{align}
    \s_n^2 = a \left[ \left( 1 - \frac{2}{\pi} \right) \sum_{j=1}^n \frac{\l_j^2}{g(j/n)} - \frac{1}{\pi n} \left( \sum_{j=1}^n \frac{\l_j}{\sqrt{g(j/n)}} \right)^2 \right].
\end{align}
For $\l = n^{-1/2}(1,\dots,1)$ both sums are Riemann sums, and since $g$ is continuous, we have
\begin{align}
    \lim_{n \to \infty} \s_n^2 = a \left[ \left( 1 - \frac{2}{\pi} \right) \int_0^1 \frac{dx}{g(x)} - \frac{1}{\pi} \left( \int_0^1 \frac{dx}{\sqrt{g(x)}} \right)^2 \right]>0. 
\end{align}
\end{enumerate}
\end{examples}
Since all assumptions are imposed coordinatewise, up to uniformity, the examples can be combined. Any mixture of the above examples, each taken with its proportion of coordinate types, in which each type occurs for a positive proportion of the coordinates is again covered, with $\b_0$ the smallest of the individual values, provided that this $\b_0$ is positive whenever some group of the coordinates needs $\b_0 > 0$ for integrability.

   \subsection*{Related results}
   Random variables under various constraints have been studied in numerous works in the literature. A classical reference to the $l_2$-sphere, the simplex and discrete simplices is \cite{DiaFr}, where finite or ``low-dimensional'' marginals are studied. This was based on normalizing i.i.d. random vectors to the set under investigation. Representation of this kind for $l_p$ balls \cite{SchechtmannZinn} led to a probabilistic approach to their geometry \cite{BartheGuedon}. Large deviations on high-dimensional $l_p$-balls were studied in \cite{GKR} while recent years have seen a strong increase in various probabilistic investigations of $l_p$-balls and related objects driven by Kabluchko, Prochno, Th\"ale \textit{et al}. Examples include the simplex as a section of the $l_1$-ball \cite{Bacietal}, Lorentz \cite{kabluchko2024probabilistic}, Schatten classes in \cite{kabluchko2020sanov}, see \cite{PTT,PPZ,ars} for more recent probabilistic, and \cite{prochno2019geometry} for a nice summary of past results. Closest to our setting geometrically are Orlicz balls $\{ \sum_j V(x_j) \leq n\}$ for a single convex function $V$, whose volumetric properties have been studied via maximum entropy and conditioning arguments \cite{KabluchkoProchno,BartheWolff}. Our surfaces are separable in the same way but we allow for non-convex and coordinate dependent $V_j$, and we study fluctuations of linear statistics under the microcanonical measure.

   From a probabilistic point of view, Theorem \ref{T1} is a conditional central limit theorem. By Remark \ref{remark_measures}, $P_n$ is the law of independent random variables conditioned on the value of a sum of functions of them. Results of this type go back to Bartlett \cite{Bartlett}, Steck \cite{Steck} and Holst \cite{Holst}, who showed that a sum conditioned on another sum is asymptotically normal, with variance given by the residual variance after regression on the conditioning sum, i.e.~by the Schur complement as in \eqref{variance}. Janson \cite{Janson} proved moment convergence in such theorems. Large and moderate deviation principles for such conditioned sums of discrete random variables were obtained in \cite{Gamboa}, also utilizing Bartlett's formula and involving the same residual variance as in \eqref{variance}. These works \cite{Holst,Janson,Gamboa} each assume lattice random variables and typically assume identical distributions. Our contribution in this direction is threefold. First, we treat non-identically distributed coordinates and weights $\l$ under assumptions that are uniform in $j$ and $n$. Second, we condition on a continuous level set for which the relevant characteristic functions decay only polynomially, controlling this decay is the main technical task of Section 3. Third, we obtain the soft constraint analog with its interpolating variance \eqref{variance2}. 
   
   In statistical physics our investigation relates to the classical question of equivalence between microcanonical and canonical ensemble as the particle number tends to infinity. While fluctuations as studied in our CLTs are not expected to be equivalent - and our theorems show indeed that they are not - they do appear as a first order correction to the equivalent expectations. A very nice paper here is \cite{CancriniOlla}, showing in a non-technical way how basic statistical mechanics principles like expansions of the partition function in the inverse temperature lead to equivalence of ensembles and accessing the variance of certain statistics. We remark that our framework does not fully coincide with a classical statistical mechanics setup in phase space. In particular, the typical kinetic energy in the momentum space introduces a Gaussian mollifier that simplifies using characteristic functions immensely as it gives rapid decay immediately \cite{CancriniOlla}. In our paper controlling the decay of characteristic functions will need a separate and lengthy technical discussion. Our focus on asymmetric surfaces also departs from the classical statistical physics framework of identical particles.

   Finally, our work fits into the broader theme of asymptotic normality of projections of high dimensional random vectors. A celebrated result is Klartag's CLT for convex bodies \cite{Klartag}, in which he shows that projections of random vectors with isotropic log-concave distributions in a high-dimensional convex body are close to a standard normal distribution for almost all directions of projection. This has been extended in \cite{KlartagBerryEsseen,Klartagrate} to Berry-Esseen bounds and convergence rates, see also \cite{Meckes, Meckes2}, and \cite{Sudakov,vonWeizsacker,Bobkov} for earlier results. In contrast, our measures live on hypersurfaces that need not be convex, need not carry log-concave measures, and may have disconnected support. Also, our results hold for every direction $\l$ with $||\l||_{\max} \to 0$, not only for almost all directions.
   
\vspace{1em}
    The rest of the paper is organized as follows. Section 2 starts the analysis of $P_n$, beginning by showing that $P_{n,\SC}$ converges weakly to $P_n$ as $\g \to \infty$, thereby providing us with a full dimensional approximation of our measure $P_n$. We also provide a very useful representation of the characteristic functions of $P_n$ and $P_{n,\SC}$  which allows us to examine the random variables under our new independent measure $P_{n,\b}$, and then show the existence of the unique $\b$ as required in Theorems $\ref{T1}$ and $\ref{T2}$. In Section 3 we provide supplementary results to be used in the main proofs of Theorems $\ref{T1}$ and $\ref{T2}$, which are done in Section 4. The most important of these supplementary results, Proposition $\ref{lemma_CF_unif}$, gives us a key bound for the characteristic function of $P_{n,\b}$.

	\section{A Representation Formula for $P_n$} 

    The main difficulties in analyzing $P_n$ of $\eqref{P_n}$, is the lower dimensionality of the support of $P_n$, making the measure seem inaccessible, and that the random variables $X_1,\dots,X_n$ are dependent, which prevents using many of the classical tools of probability theory. The main result of this section will be a formula that does not only turn the lower-dimensional problem into a full-dimensional problem, but that also allows us to use independent random variables, albeit under a different distribution. The idea is that formally the constraint \eqref{constraint} can be understood as a Dirac ``delta function'' 
	\begin{align} 
		\d\lb \sum_{j=1}^n V_j(X_j)-an\rb. 
	\end{align} 
	This is made rigorous by approximating this delta function using the Gaussian term
	\begin{align} 
		\sqrt{\frac\g{\pi n}} e^{-\g\lb \frac1{\sqrt n}\sum_{j=1}^n V_j(x_j)-a\sqrt n\rb^2}
	\end{align} 
	and sending $\g\to\infty$. This provides the motivation for the definition of $P_{n,\SC}$ in \eqref{P_n,SC}. 
    It will be convenient to rewrite the density of $P_{n,\SC}$ in order to make the integrability of the density clearer. For any $\b \in \R$ we can write
	\begin{align} 
		-\g\lb \frac1{\sqrt n}\sum_{j=1}^n V_j(x_j)-a\sqrt n\rb^2=&-\g\lb \frac1{\sqrt n}\sum_{j=1}^n V_j(x_j)-(a+\frac\b{2\g})\sqrt n\rb^2-\b\sum_{j=1}^nV_j(x_j)\\ 
		&+\frac{\b^2n}{4\g}+\b an.\label{manipulation} 
	\end{align} 
	As the last two terms are constant, they can be absorbed into the normalizing constant. We note in passing that taking $\b = \b_0$ and bounding the first term on the right by zero, ensures finiteness of $Z_{n,\SC}$ for all fixed $n$ by Assumption \ref{ass:int}. In particular the $f_j$ need not be integrable themselves, for instance taking $f_j =1$ is allowed as long as Assumption \ref{ass:int} holds for some $\b_0 \geq 0$.

    We now make these ideas rigorous in three propositions, whose proofs are deferred to the end of this section. The first concerns the convergence of $P_{n,\SC}$ to $P_n$.

    \begin{prop}\label{prop:convergence} 
        Let Assumptions \ref{ASS} and \ref{MASS} hold and recall the definition of $G$ and $p$ from Assumption \ref{ass:p}. For $n$ sufficiently large such that $|G| > p+1$, $P_{n,\SC}$ converges weakly to $P_n$ as $\g\to\infty$. 
	\end{prop} 

    The following proposition provides a characteristic function representation for $P_n$ and $P_{n,\SC}$ in terms of $P_{n,\b}$ of \eqref{P_n,b}. 
    
	\begin{prop}\label{prop_CF}
     Let Assumptions $\ref{ASS}$ and $\ref{MASS}$ hold.
		\begin{enumerate} 
        \item For any $\l\in\R^n$, $\b \geq \b_0$, n sufficiently large such that $|G| >p+1$, and all $t\in\R$ we have
			\begin{align}	
				&\E_{{n}}\lb e^{it\sum_{j=1}^n\l_j(X_j-\E_{n,\b}X_j)}\rb \\
				&= \frac{\int_{\R}\E_{{n,\b}}\lb e^{is\lb \frac1{\sqrt n}\sum_{j=1}^n V_j(X_j)-a\sqrt n\rb + it\sum_{j=1}^n\l_j(X_j-\E_{n,\b}X_j)}\rb ds}{\int_{\R}\E_{{n,\b}}\lb e^{is\lb \frac1{\sqrt n}\sum_{j=1}^n V_j(X_j)-a\sqrt n\rb}\rb ds}.   
                 \label{2.32}				    
			\end{align}
			\item  For any $\l\in\R^n$, $\b \geq \b_0$, $n\in\mathbb N$, and all $t\in\R$ we have
			\begin{align}
				&\E_{{n,\SC}}\lb e^{it\sum_{j=1}^n\l_j(X_j-\E_{n,\b}X_j)}\rb \\
				&= \frac{\int_{\R}\E_{{n,\b}}\lb e^{is\lb \frac1{\sqrt n}\sum_{j=1}^n V_j(X_j)-\left(a+\frac{\b}{2\g}\right)\sqrt n\rb + it\sum_{j=1}^n\l_j(X_j-\E_{n,\b}X_j)}\rb e^{-\frac{1}{4\g}s^2} ds}{\int_{\R}\E_{{n,\b}}\lb e^{is\lb \frac1{\sqrt n}\sum_{j=1}^n V_j(X_j)-\left(a+\frac{\b}{2\g}\right)\sqrt n\rb}\rb e^{-\frac{1}{4\g}s^2} ds}.  \label{2.31}
			\end{align} 
		\end{enumerate} 
	\end{prop} 
\vspace{1em}
\begin{remark}\label{remark_Bartlett}
    Proposition \ref{prop_CF} allows us to work in a set of full dimensionality and under a product measure $P_{n,\b}$, albeit at the expense of an additional integral. With this representation classical tools of independence become available. We also see very clearly that the r.h.s.~of, for example, part \eqref{2.32} of Proposition \ref{prop_CF}, involves two sums at different scales, the weighted $\sum_{j=1}^n\l_j(X_j-\E_{n,\b}X_j)$ with weights $\l_j$, and the classically weighted $\frac1{\sqrt n}\sum_{j=1}^n V_j(X_j)-a\sqrt n$. This interplay of these two statistics leads to the shape of the variance in Theorem \ref{T1}. Formula \eqref{2.32} is a variant of a formula due to Bartlett \cite{Bartlett}, derived for conditional characteristic functions of random variables, which has been used similarly in \cite{Gamboa} in a discrete setting. 
\end{remark}

	We will now turn to the choice of $\b$. In the case of the hard constraint, the optimal $\b$ will be the one such that under $P_{n,\b}$ the constraint is fulfilled in expectation, i.e.~$\E_{n,\b}(\sum_{j=1}^n V_j(X_j))=an$. For the soft constraint case things are slightly different. In line with what has been discussed in \eqref{manipulation} the aim is now to choose $\b\in\R$ such that under $P_{n,\b}$ the modified constraint 
    \begin{align}
	    \frac1{\sqrt n}\sum_{j=1}^n V_j(X_j)=(a+\frac\b{2\g})\sqrt n\label{SC_beta}
	\end{align} 
    is satisfied in expectation. The following proposition provides the existence and uniqueness of such a $\b$, and also ensures that this $\b$ is bounded uniformly in $n$. 
    
	\begin{prop} \label{lemma_beta} 
		Let Assumptions $\ref{ASS}$  and $\ref{MASS}$ hold. Then there is a unique $\b\ge\b_0$ such that $\E_{{n,\b}}\left(\sum_{j=1}^n V_j(X_j)\right)=na$. Similarly, for each $\g>0$ there is a unique $\b=\b(\g)\geq \b_0$ such that $\E_{{n,\b}}\left(\sum_{j=1}^n V_j(X_j)\right)=\left(a+\frac{\b}{2\g}\right)n$. Moreover, in both cases there is a constant $\b^*< \infty$, independent of $n$, with 
		\begin{align}
			\b_0 \leq \b \leq \b^*\quad  \text{for all }n.
			\end{align}
	\end{prop} 
\vspace{1em}
In order to provide the proofs of the preceding propositions we first provide a lemma which gives an estimate on the one dimensional oscillatory integrals that will be used repeatedly throughout the paper. For a single coordinate $j$, 
\begin{align}
    \int_{\R} e^{irV_j(x)} f_j(x) e^{-\b V_j(x) + i\tau x} dx 
\end{align}
is, up to normalization, the joint characteristic function of $V_j(X_j)$ and $X_j$ under $P_{n,\b}$, evaluated at $(r,\t)$. Sufficient decay of this integral will be needed in Proposition \ref{prop:convergence} to show that our measures are well defined, in Proposition \ref{prop_CF} to interchange limit $\g \to \infty$ and integration and in Proposition \ref{lemma_CF_unif} to deal with the characteristic function at large frequencies. The following lemma provides the decay uniform in $j$ and $n$.
   \begin{lemma} \label{lem_vdc}
        For all $j \in G$ of Assumption \ref{ass:p}, and $\b \geq \b_0$ of Assumption \ref{ass:BM}, $\tau \in \R$, and $|r| \geq 1$,
        \begin{align}
            \lvv \int_{\R} e^{irV_j(x)} f_j(x) e^{-\b V_j(x) + i\tau x} dx \rvv \leq C_0\lb 1+\b+|\t|\rb|r|^{-1/p},
        \end{align}
        with $C_0$ depending only on constants in Assumptions 1 and 2, in particular $C_0$ is uniform in $j$ and $n$.
\end{lemma}
\begin{proof}[Proof of Lemma \ref{lem_vdc}]
Fix $j \in G$, $\b \geq \b_0$, $\t$, and $r$ and set 
\begin{align}
    \psi_j(x) := f_j(x)e^{-\b V_j(x) + i \t x}.
\end{align}
Since $V_j \geq 0$ and $\b \geq \b_0$, we have $|\psi_j(x)| \leq M$ by Assumption \ref{ass:regularity} and $|\psi_j(x)|$ integrable by Assumption \ref{ass:int}. So the integral converges absolutely. Recalling $D_j$ and $H_j$ as the set of points of discontinuity of $f_j$ and non-differentiability of $V_j$ respectively as in Assumption \ref{ass:regularity}, on $\R \setminus (D_j \cup H_j)$ the function $\psi_j$ is continuously differentiable, and we can write
\begin{align} 
 &\int_{\R\setminus (D_j \cup H_j)} |\psi_j'(x)| dx \leq \int_{\R\setminus (D_j \cup H_j)} \lb |f_j'(x)| + \b |V'_j(x)f_j(x)| +  |\t| |f_j(x)| \rb e^{-\b V_j(x)} dx\\
 &\leq K(1+ \b +|\tau| ),
\end{align}
using $e^{-\b V_j} \leq e^{-\b_0 V_j}$ and Assumption \ref{ass:int}.
Let $I_{j,1},\dots,I_{j,A_j}$ be the intervals of Assumption \ref{ass:p} with $A_j \leq T$, and define the disjoint intervals $ \tilde I_{j,\tilde k}:= I_{j,k} \setminus \bigcup_{l<k} I_{j,l}$. Each of these is a union of at most $k \leq T$ intervals, so there are at most $T^2$ intervals in total, and each interval is a subset of some $I_{j,k}$. Cutting each interval $\tilde I_{j,\tilde k}$ at the points of $(D_j \cup H_j)$ gives disjoint open intervals $J_{j,1},\dots,J_{j,N}$, with $N$ uniformly bounded in $j$. On each $J_{j,i}$ $i \in \{1,\dots,N\}$, $\psi_j$ is continuously differentiable, and since each $J_{j,i} \subseteq \tilde I_{j,\tilde k} \subseteq I_{j,k}$ for some $k$, Assumption \ref{ass:p} ensures there is an order $k_i \in \{1,\dots,p\}$ with $|V_j^{(k_i)}(x)| \geq \e$ on $J_{j,i}$. We also note that if $J_{j,i}$ is unbounded, then so is the corresponding $I_{j,k}$, and hence $k_i = 1$ by Assumption \ref{ass:p}. For intervals with $k_i \geq 2$, a standard application of the van der Corput lemma as in \cite[Chapter VIII, Proposition 2 and Corollary]{Stein} gives us that 
\begin{align}
     \lvv \int_{J_{j,i}} e^{irV_j(x)} f_j(x) e^{-\b V_j(x) + i\tau x} dx \rvv \leq C\lb 1+\b+|\t|\rb|r|^{-1/k_i},
\end{align}
for some $C>0$ uniform in $j$.
For intervals for which $k_i =1 $ we have $|V_j'(x)|>\e$ and we can use integration by parts as well as Assumption \ref{ass:int} to give us
\begin{align}
      \lvv \int_{J_{j,i}} e^{irV_j(x)} f_j(x) e^{-\b V_j(x) + i\tau x} dx \rvv \leq \lb \frac{2M + K(1+ \b +|\tau| )}{\e} + K\e^{-2}\rb |r|^{-1}.
\end{align}
Finally, since $|r|\geq 1$ and $k_i \leq p$, we have $|r|^{-1/k_i} \leq |r| ^{-1/p}$. Summing over the intervals, with the number of intervals being uniformly bounded in $j$, completes the proof of the lemma.
\end{proof}
\vspace{1em}

Equipped with Lemma \ref{lem_vdc} we can now present the proofs of Propositions \ref{prop:convergence}, \ref{prop_CF} and \ref{lemma_beta}. We begin with the proof of Proposition \ref{prop:convergence} showing that $P_{n,\SC} \to P_n$ weakly as $\g \to \infty$. 

	\begin{proof}[Proof of Proposition \ref{prop:convergence}]
		We will show that the distribution function of $P_{n,\SC}$ converges pointwise to that of $P_n$.
        \newline \indent 
		 To begin let $E = \bigtimes_{j=1}^n ( -\infty , b_j],$ $ b_j\in\R$. By the elementary identity \eqref{manipulation} above with $\b = \b_0$ we can write
		\begin{align} 
			&P_{n,\SC}(E)=\frac{e^{\b_0 an} \sqrt{\frac\g{\pi n}}\int_{E}e^{-\g \lb\frac{1}{\sqrt n} \sum_{j=1}^n V_j(x_j) - \lb a + \frac{\b_0}{2\g}\rb\sqrt n \rb^2}\prod_{j=1}^n f_j(x_j)e^{-\b_0 V_j(x_j)}dx} 
			{e^{\b_0 an}\sqrt{\frac\g{\pi n}}\int_{\R^n}e^{-\g \lb\frac{1}{\sqrt n} \sum_{j=1}^n V_j(x_j) - \lb a + \frac{\b_0}{2\g}\rb\sqrt n \rb^2}\prod_{j=1}^n f_j(x_j)e^{-\b_0 V_j(x_j)}dx},
            \label{Laplace1} 
		\end{align} 
        and by Assumption \ref{ass:int} both integrals are finite.
        \newline \indent
        The proof consists of two main steps. First, we evaluate the limit as $\g \to \infty$ of the numerator of \eqref{Laplace1}, leading to the expression \eqref{distr_func_lhs}. Then we will use the area formula to relate this limit with $Z_n P_n(E)$. The denominator can be treated the same as the numerator with $E = \R^n$.
        To begin we wish to make a change of variables in the $n$-th coordinate and to this end we assume without loss of generality that $n \in G$ as in Assumption \ref{ass:p}. Examining the numerator of \eqref{Laplace1} and integrating over $x_n$ first, the domain $(-\infty,b_n]$ of the one-dimensional integral can be split into $m$ open intervals $E_k$, on which $V_n$ is continuously differentiable with $V'_n \neq 0$, by cutting at the points of $H_n$, as in Assumption \ref{ass:regularity}, $b_n$ and the stationary points of $V_n$, of which there are finitely many on the support of $f_n$ by Assumption \ref{ass:p}, and elsewhere the integrand is zero. For each interval $E_k$ the function $x_n\mapsto V_n(x_n)$ is one-to-one, allowing us to apply the change of variables $\psi_k: E_k\to\R$,		
\begin{align} 
	y_n:=\psi_k(x_n):=\sum_{j=1}^n  V_j(x_j),
		\end{align} 
        to the numerator of \eqref{Laplace1}. Defining
        \begin{align}
			q(u): = \sum_{k=1}^m f_n(V_n^{-1}(u))e^{-\b_0u} \dopp1_{V_n(E_k)}(u)\lv V_n'(V_n^{-1}(u))\rv^{-1}\label{indicator},
		\end{align}
        with $V_n^{-1}$ denoting the inverse of $V_n$ on $E_k$, and applying this change of variables allows us to write the numerator of \eqref{Laplace1} as
		\begin{align} 
			&e^{\b_0 an} \sqrt{\frac\g{\pi n}}\\
            &\times\int_{\R^{n-1}}\prod_{j=1}^{n-1} f_j(x_j)e^{-\b_0 V_j(x_j)}\dopp1_{(-\infty,b_j]}(x_j) \int_\R q\lb y_n - \sum_{j=1}^{n-1} V_j(x_j)\rb e^{-\g\lb\frac{y_n}{\sqrt n} - \lb a + \frac{\b_0}{2\g}\rb\sqrt n \rb^2} dy_ndx_1\dots dx_{n-1},\label{Laplace2} 
		\end{align}
        where we have used that 
        \begin{align}
            \psi_k^{-1}(y_n) = V_n^{-1} \lb y_n - \sum_{j=1}^{n-1} V_j(x_j) \rb.
        \end{align}
        We note that applying the same change of variables in the opposite direction allows us to write
        \begin{align}
            &\int_\R q(u) du =  \sum_{k=1}^m   \int_{V_n(E_k)} f_n(V_n^{-1}(u))e^{-\b_0u} \lv V_n'(V_n^{-1}(u))\rv^{-1}du
            \\
            &= \sum_{k=1}^m \int_{E_k}f_n(x_n) e^{-\b_0 V_n(x_n)} dx_n = \int_{-\infty}^{b_n} f_n(x_n) e^{-\b_0 V_n(x_n)} dx_n \leq K \label{mbleint}
        \end{align}
        with Assumption \ref{ass:int} providing the upper bound.
        \newline \indent
        Next we show that the first $n-1$ co-ordinates have a bounded density. We consider the finite measure $\s_1$, with
        \begin{align}
            \sigma_1 (A) := \int_{\R^{n-1}} \dopp1_{A} \lb \sum_{j=1}^{n-1} V_j(x_j)\rb \prod_{j=1}^{n-1} f_j(x_j) e^{-\b_0 V_j(x_j)} \dopp 1_{(-\infty,b_j]}(x_j) dx_1\dots dx_{n-1},
        \end{align}
        for all Borel sets $A$. That is $\s_1$ is the image of the finite measure $\prod_{j=1}^{n-1} f_j(x) e^{-\b_0 V_j(x)} \dopp 1_{(-\infty,b_j]}(x_j) dx$ under the map $x \mapsto \sum_{j=1}^{n-1} V_j(x_j)$. For a finite Borel measure $\s$ on $\R$, we denote its Fourier transform by
        \begin{align}
            \hat \s(r) := \int_\R e^{iry} \s(dy), \quad r \in \R.
        \end{align}
        By definition of $\s_1$ we have
        \begin{align}
            \hat \sigma_1 (r) = \prod_{j=1}^{n-1} \int_{-\infty}^{b_j} e^{ir V_j(x)} f_j(x) e^{-\b_0 V_j(x)}dx = \O\lb |r|^{-\frac{|G\setminus\{n\} | }{p}} \rb
        \end{align}
        where the bound is supplied by Lemma \ref{lem_vdc}, applied to $j \in G\setminus \{n\}$ with $\t=0$, $\b= \b_0$, and noting that for $j \not \in G\setminus\{n\}$ the factors are bounded by $K$ as in Assumption \ref{ass:int}. This is integrable since $|G\setminus\{n\}| > p$ by hypothesis, hence the measure $\s_1$ has a density $s:\R \to [0,\infty)$ which is bounded. Then for every measurable $H: \R\to\R_+$ the change of variables formula for image measures tells us that
        \begin{align}
            &\int_{\R^{n-1}} H \lb \sum_{j=1}^{n-1} V_j(x_j) \rb  \prod_{j=1}^{n-1} f_j(x_j) e^{-\b_0 V_j(x_j)} \dopp 1_{(-\infty,b_j]}(x_j) dx_1\dots dx_{n-1}\\
            &= \int_\R H(v) s(v) dv.\label{StoG}
        \end{align}
        Here $v$ represents the value of $\sum_{j=1}^{n-1} V_j(x_j)$.
        \newline \indent
        We can use \eqref{Laplace2} to rewrite the numerator of \eqref{Laplace1} as 
        \begin{align}
        &e^{\b_0 an} \sqrt{\frac\g{\pi n}}\\
        &\times \int_{\R^{n-1}}\prod_{j=1}^{n-1} f_j(x_j)e^{-\b_0 V_j(x_j)}\dopp1_{(-\infty,b_j]}(x_j) \int_\R q\lb y_n - \sum_{j=1}^{n-1} V_j(x_j)\rb e^{-\g\lb\frac{y_n}{\sqrt n} - \lb a + \frac{\b_0}{2\g}\rb\sqrt n \rb^2} dy_ndx_1\dots dx_{n-1} \\
        & = e^{\b_0 an} \\
        &\int_{\R^{n-1}} \left[\int_{\R} \sqrt{\frac \g{\pi n}} e^{-\g  \lb\frac{y_n}{\sqrt n} - \lb a + \frac{\b_0}{2\g}\rb\sqrt n \rb^2} q(y_n-\sum_{j=1}^{n-1}V_j(x_j))dy_n\right] \prod_{j=1}^{n-1} f_j(x_j) e^{-\b_0 V_j(x_j)} \dopp 1_{(-\infty,b_j]}(x_j) dx_1\dots dx_{n-1}.
        \end{align}
        We now wish to apply the change of variables formula for image measures, as in \eqref{StoG} with $H$ defined as
        \begin{align}
            H(v) : = \int_{\R} \sqrt{\frac \g{\pi n}} e^{-\g  \lb\frac{y_n}{\sqrt n} - \lb a + \frac{\b_0}{2\g}\rb\sqrt n \rb^2} q(y_n-v)dy_n,
        \end{align}
        for $v \in \R$.
        Since $q \geq 0$ is measurable, $H$ is a nonnegative measurable function and it can be seen from \eqref{mbleint} that it is finite. 
        Now utilizing the change of variables formula for image measures, as in \eqref{StoG}, with this $H$, the numerator of \eqref{Laplace1} equals
       \begin{align}
          &e^{\b_0 an} \int_{\R^{n-1}} \left[\int_{\R} \sqrt{\frac \g{\pi n}} e^{-\g  \lb\frac{y_n}{\sqrt n} - \lb a + \frac{\b_0}{2\g}\rb\sqrt n \rb^2} q(y_n-v)dy_n\right] \prod_{j=1}^{n-1} f_j(x_j) e^{-\b_0 V_j(x_j)} \dopp 1_{(-\infty,b_j]}(x_j) dx_1\dots dx_{n-1}\\
          &= e^{\b_0 an} \int_{\R} s(v) \int_{\R} \sqrt{\frac \g{\pi n}}e^{-\g \lb\frac{y_n}{\sqrt n} - \lb a + \frac{\b_0}{2\g}\rb\sqrt n \rb^2}q(y_n-v) dy_n dv.
       \end{align}
       The integrand is nonnegative, so by Tonelli's theorem we may interchange order of integration. The numerator of \eqref{Laplace1} therefore equals
       \begin{align}
           e^{\b_0 an} \int_\R \sqrt{\frac{\g}{\pi n}} e^{-\g \lb\frac{y_n}{\sqrt n} - \lb a + \frac{\b_0}{2\g}\rb\sqrt n \rb^2} \int_\R s(v) q(y_n-v) dv dy_n. \label{Laplace3}
       \end{align}
       \newline \indent
       We define the function 
       \begin{align}
           F(y) : = \int_\R s(v)q(y-v)dv,
       \end{align}
       such that \eqref{Laplace3} can be written as 
       \begin{align}
              e^{\b_0 an} \int_\R \sqrt{\frac{\g}{\pi n}} e^{-\g \lb\frac{y}{\sqrt n} - \lb a + \frac{\b_0}{2\g}\rb\sqrt n \rb^2} F(y) dy.
       \end{align}
       The function $F$ can be understood as the unnormalized density of the constraint level $\sum_{j=1}^n V_j(x_j)$, and is bounded and continuous even though $q$ need not be. Indeed $0\leq F \leq K ||s||_{\max}$, from \eqref{mbleint}, and for all $y,y' \in \R$,
       \begin{align}
           |F(y)-F(y')| \leq ||s||_{\max} \int_\R |q(u+y-y')-q(u)|du,
       \end{align}
       which tends to $0$ as $y - y' \to 0$, by the continuity of translation in $L^1(\R)$.
       We note that
       \begin{align}
           \sqrt{\frac{\g}{\pi n}}\exp\lb-\g \lb\frac{y}{\sqrt n} - \lb a + \frac{\b_0}{2\g}\rb\sqrt n \rb^2\rb
       \end{align}
       is the density of a normal distribution with mean $m_\g:= (a+\b_0/(2\g))n$ and variance $\s^2_\g := n/(2\g)$. Substituting $y = m_\g + \s_\g z$, the numerator of \eqref{Laplace1} equals
       \begin{align}
           e^{\b_0 an} \int_\R \frac{e^{-\frac{z^2}{2}}}{\sqrt{2\pi}} F(m_\g + \sigma_\g z) dz.
       \end{align}
       As $\g \to \infty$ we have $m_\g +\s_\g z \to an$ for every $z \in \R$, so $F(m_\g + \s_\g z) \to F(an)$ by continuity of $F$. Since the integrand is bounded by $K||s||_{\max} e^{-z^2 /2}/\sqrt{2\pi}$, uniformly in $\g$, Lebesgue's dominated convergence theorem lets us write
       \begin{align}
           &\lim_{\g \to \infty} e^{\b_0 an}  \int_\R e^{-\frac{z^2}{2}}\lb\sqrt{2\pi}\rb^{-1} F(m_\g + \sigma_\g z) dz = e^{\b_0 an} F(an)= e^{\b_0 an} \int_\R s(v) q(an -v ) dv
          \\  &= e^{\b_0 an} \int_{\R^{n-1}} q \lb an - \sum_{j=1}^{n-1} V_j(x_j) \rb  \prod_{j=1}^{n-1} f_j(x_j) e^{-\b_0 V_j(x_j)} \dopp 1_{(-\infty,b_j]}(x_j) dx_1\dots dx_{n-1}
       \end{align}
       where we again use \eqref{StoG}. 
       \newline \indent 
       Defining
       \begin{align}
           c(x_1,\dots,x_{n-1}) := \sum_{j=1}^{n-1} V_j(x_j),
       \end{align} 
       and recalling the definition of $q$ in \eqref{indicator} we can write
		\begin{align}
			&q\lb an- c\rb = \sum_{k=1}^m f_n\lb V_n^{-1}\lb an- c\rb \rb e^{-\b_0\lb an- c\rb} \dopp1_{V_n(E_k)}\lb an-c\rb \lv V_n'\lb V_n^{-1}\lb an- c\rb \rb\rv^{-1}.
		\end{align}
        Since by definition $\psi_k(x_n) = c+V_n(x_n)$ on $E_k$, we have $an -c \in V_n(E_k)$ if and only if $an \in \psi_k(E_k)$, and then $V_n^{-1}(an-c) = \psi_k^{-1}(an)$. Hence, we can write
        \begin{align}
            q\lb an - \sum_{j=1}^{n-1} V_j(x_j) \rb =  \sum_{k=1}^m f_n(\psi_k^{-1}(an))\dopp1_{\psi_k(E_k)}(an)\lv V_n'(\psi_k^{-1}(an))\rv^{-1} e^{-\b_0 (an -c)} .
        \end{align}
        We observe that
		\begin{align}
		\prod_{j=1}^{n-1}\dopp1_{(-\infty,b_j]}(x_j)\cdot	\dopp1_{\psi_k(E_k)}(an)=\dopp1_{\textup{Proj}((\bigtimes_{j=1}^{n-1}(-\infty,b_j]\times E_k)\cap S_{a,V})}(x_1,\dots,x_{n-1}),
		\end{align} 
	where Proj denotes the projection from $\R^n$ to $\R^{n-1}$ onto the first $n-1$ coordinates. We thus obtain that the limit $\g\to\infty$ of \eqref{Laplace2} is 
	\begin{align}
		 \sum_{k=1}^m &\int_{\R^{n-1}}\prod_{j=1}^{n-1} f_j(x_j)\\
		\times& f_n(\psi_k^{-1}(an))\lv V_n'(\psi_k^{-1}(an))\rv^{-1}\dopp1_{\textup{Proj}((\bigtimes_{j=1}^{n-1}(-\infty,b_j]\times E_k)\cap S_{a,V})}(x_1,\dots,x_{n-1})dx_1\dots dx_{n-1}\label{distr_func_lhs},
	\end{align}
    where we used that, with $c = \sum_{j=1}^{n-1} V_j(x_j) $ the exponential factors combine to
    \begin{align}
        e^{\b_0 an} e^{-\b_0 (an-c)} \prod_{j=1}^{n-1} e^{-\b_0 V_j(x_j)}=1.
    \end{align}
    \newline \indent
    We now attempt to understand \eqref{distr_func_lhs} as 
	\begin{align}
		Z_n P_n(E)=\int_{\R^{n}}\prod_{j=1}^{n} f_j(x_j)\dopp1_{E\cap S_{a,V}}(x) \frac{1}{\|\nabla V(x)\|}dH^{n-1}(x),
	\end{align}
	where we recall $V:\R^n\to\R$, $V(x):=\sum_{j=1}^n V_j(x_j)$.
To do that, we employ the area formula, expressing an $n$-dimensional integral w.r.t.~the $n-1$-dimensional Hausdorff measure as an $n-1$-dimensional integral w.r.t.~the $n-1$-dimensional Lebesgue measure. Before we do this we recall the definition of $H_j$ as in Assumption \ref{ass:regularity} and define the set
\begin{align}
    &S : =\\
    &\{ x \in S_{a,V} \cap \supp (\prod_{j=1}^n f_j): V'_n(x_n) = 0 \text{ or } x_n \in D_n \cup \{b_n\}\text{ or } x_j \in H_j \text{ for some }j \in \{1,\dots, n\}\}. \label{S1}
\end{align}
Since $n \in G$, $V_n$ has finitely many stationary points on $\supp(f_n)$, and $|H_j| +|D_j| \leq T$ by Assumption \ref{ass:regularity}. Hence $S$ is contained in the union of slices $\{x_i = z\} \cap S_{a,V}$ for finitely many pairs $(i,z)$, with $i \in \{1,\dots,n\}$. On a given slice the remaining coordinates must satisfy $\sum_{j \neq i} V_j(x_j) = an - V_i(z)$, since $x \in S_{a,V}$. Choosing $j \in G\setminus\{i\}$, for fixed values of the other coordinates this equation has only finitely many solutions for $x_j \in \supp(f_j)$, by Assumption \ref{ass:p}. By Fubini's theorem each slice is therefore a Lebesgue null set in the hyperplane $\{x_i =z \}$, and hence $H^{n-1}(S) = 0$.

For $k \in \{1,\dots,m\}$ we define the parameterization
\begin{align}
	\Phi_k(x_1,\dots,x_{n-1}):=\lb x_1,\dots,x_{n-1},V_n^{-1}(an-\sum_{j=1}^{n-1}V_j(x_j))\rb=\lb x_1,\dots,x_{n-1},\psi_k^{-1}(an)\rb,
\end{align}
with $V_n^{-1}$ the inverse of $V_n$ on $E_k$ as mentioned earlier. This parameterization is one-to-one on $\textup{Proj}((\bigtimes_{j=1}^{n-1}(-\infty,b_j]\times E_k)\cap S_{a,V})$, discarding any empty sets. It is continuously differentiable off the Lebesgue null set $S':=\{x \in \R^{n-1}:x_j \in H_j \text{ for some } j < n\}$. In the area formula applied below neither side is changed by removing $S'$ from the integrals. In particular, the integral over $\R^{n-1}$ is unchanged since $S'$ is Lebesgue null, and the integral over $S_{a,V}$ is unchanged because $\Phi_k(S') \subset S$ and $H^{n-1}(S) = 0$. The Jacobian matrix has the form
\begin{align}
	D\Phi_k(x)=\begin{pmatrix}
		I_{n-1}\\
		v^T
	\end{pmatrix},
\end{align}
where $I_{n-1}$ is the $n-1$-dimensional identity matrix and $v$ is the $n-1$-dimensional column vector with 
\begin{align}
	v_j:=-\frac{V_j'(x_j)}{V_n'(V_n^{-1}(an-\sum_{j=1}^{n-1}V_j(x_j)))},\quad j=1,\dots,n-1.
\end{align}
It is then straightforward to check that the Jacobian for the area formula is
\begin{align}
	J\Phi_k(x)&=\sqrt{\det(D\Phi_k(x)^TD\Phi_k(x))}=\sqrt{\det(I_{n-1}+vv^T)}=\sqrt{1+\|v\|^2}\\
	&=\sqrt{1+\frac{\sum_{j=1}^{n-1}V_j'(x_j)^2}{V_n'(V_n^{-1}(an-\sum_{j=1}^{n-1}V_j(x_j)))^2}}=\sqrt{\frac{\sum_{j=1}^{n}V_j'(x_j)^2}{V_n'(x_n)^2}}=\frac{\|\nabla V(x)\|}{\lv V_n'(x_n)\rv},
\end{align}
where we define $x_n:=V_n^{-1}(an-\sum_{j=1}^{n-1}V_j(x_j))$. We now obtain with the area formula  \cite[Theorem 3.2.3]{Federer} and using that $S$ from \eqref{S1} has $H^{n-1}$-measure 0, that
\begin{align}
	&Z_n P_n(E)=\int_{\R^{n}}\prod_{j=1}^{n} f_j(x_j)\dopp1_{E\cap S_{a,V}}(x) \frac{1}{\|\nabla V(x)\|}dH^{n-1}(x)\\
	&=\sum_{k=1}^m\int_{\R^n}\prod_{j=1}^{n} f_j(x_j)\dopp1_{(\bigtimes_{j=1}^{n-1}(-\infty,b_j]\times E_k)\cap S_{a,V}}(x) \frac{1}{\|\nabla V(x)\|}dH^{n-1}(x)\\
		&=\sum_{k=1}^m \int_{\R^{n-1}}\prod_{j=1}^{n-1} f_j(x_j)\frac{1}{\|\nabla V(x)\|}J\Phi_k(x)\\
	\times& f_n(\psi_k^{-1}(an))\dopp1_{\textup{Proj}((\bigtimes_{j=1}^{n-1}(-\infty,b_j]\times E_k)\cap S_{a,V})}(x_1,\dots,x_{n-1})dx_1\dots dx_{n-1},
\end{align} 
which is \eqref{distr_func_lhs}. The denominator of \eqref{Laplace1} is the numerator with $E$ replaced by $\R^n$, that is, with $b_j = \infty$ for all $j$. Running the same argument, with $F_\R$ denoting the function $F$ obtained for this choice, shows that the denominator converges to $e^{\b_0 an}F_\R(an) = Z_n$, which is finite since $F_\R$ is bounded and positive by assumption. Hence $P_{n,\SC}(E) \to Z_n P_n(E)/Z_n = P_n(E)$, which completes the proof.
	\end{proof} 
\vspace{1em}	
We now give the proof of Proposition \ref{prop_CF}. We will first prove the soft constraint case and the hard constraint case will follow from taking $\g \to \infty$ noting that we have just proved that as $\g \to \infty$ $P_{n,\SC} \to P_n$ weakly, and hence the characteristic function converges pointwise.

	\begin{proof}[Proof of Proposition \ref{prop_CF}] 
		We start with proving $\eqref{2.31}$, the soft constraint case. Using the density of $P_{n,\SC}$ alongside $\eqref{manipulation}$ we have with the shorter notation $\sum_{j=1}^n\l_j(X_j-\E_{n,\b}X_j)=\langle \l , X - \E_{n,\b}X\rangle$, that
		\begin{align} 
			&\E_{{n,\SC}}(e^{it\langle \l , X - \E_{n,\b}X\rangle})\\
            &= \int_{\R^n} \frac{e^{it\langle \l , x - \E_{n,\b}X\rangle}}{Z_{n,\SC}}e^{-\g\left(\frac{1}{\sqrt{n}}\sum_{j=1}^{n}V_j(x_j) - (a +\frac{\b}{2\g})\sqrt{n}\right)^2+\frac{\b^2n}{4\g}+\b an} \prod_{j=1}^{n} e^{-\b V_j(x_j)}f_j(x_j)dx, \label{ESC}
		\end{align} 
		and we can apply the Gaussian Fourier transform identity
		\begin{align} 
			e^{-\g z^2}=\frac{1}{\sqrt{4\pi\g}}\int_{\R}e^{isz}e^{-\frac{1}{4\g}s^2}ds \label{LNRSTN} 
		\end{align} 
		to rewrite $\eqref{ESC}$ as
		\begin{align} 
			& \frac{e^{\frac{\b^2n}{4\g}+\b an}}{\sqrt{4\pi\g}}  \frac{1}{Z_{n,\SC}} \\
            &\times\int_{\R} \int_{\R^n}  e^{it\langle \l , x - \E_{n,\b}X\rangle + is\left(\frac{1}{\sqrt{n}}\sum_{j=1}^{n}V_j(x_j) - \left(a +\frac{\b}{2\g}\right)\sqrt{n}\right)-\frac{1}{4\g}s^2}\prod_{j=1}^{n} e^{-\b V_j(x_j)}f_j(x_j) dxds .  \label{ESC1}
		\end{align} 
		We now observe that 
		\begin{align} 
			\frac{Z_{n,\SC}}{Z_{n,\b}} = \E_{{n,\b}}\left(e^{-\g\left(\frac{1}{\sqrt{n}}\sum_{j=1}^{n}V_j(X_j) - (a +\frac{\b}{2\g})\sqrt{n}\right)^2+\frac{\b^2n}{4\g}+\b an}\right) ,
		\end{align} 
		and multiply $\eqref{ESC1}$ by ${Z_{n,\b}}/{Z_{n,\b}}$ to yield 
		\begin{align} 
			&\E_{{n,\SC}}(e^{it\langle \l , X - \E_{n,\b}X\rangle})  \\
			&=\frac{e^{\frac{\b^2n}{4\g}+\b an}\frac{1}{\sqrt{4\pi\g}} \int_{\R} \E_{n,\b}\left(e^{it\langle \l , X - \E_{n,\b}X\rangle +is\left(\frac{1}{\sqrt{n}}\sum_{j=1}^{n}V_j(X_j) - \left(a +\frac{\b}{2\g}\right)\sqrt{n}\right)-\frac{1}{4\g}s^2}\right) ds }{\E_{{n,\b}}\left(e^{-\g\left(\frac{1}{\sqrt{n}}\sum_{j=1}^{n}V_j(X_j) - \lb a +\frac{\b}{2\g}\rb \sqrt{n}\right)^2+\frac{\b^2n}{4\g}+\b an}\right)} \label{ESC2}.
		\end{align} 
		Now applying $\eqref{LNRSTN}$ to the denominator of $\eqref{ESC2}$ yields 
		\begin{align} 
			\E_{{n,\SC}}(e^{it\langle \l , X - \E_{n,\b}X\rangle}) = \frac{\int_{\R}\E_{{n,\b}}\lb e^{is\lb \frac1{\sqrt n}\sum_{j=1}^n V_j(X_j)-\left(a+\frac{\b}{2\g}\right)\sqrt n\rb + it\langle \l,X-\E_{n,\b} X \rangle}\rb e^{-\frac{1}{4\g}s^2} ds}{\int_{\R}\E_{{n,\b}}\lb e^{is\lb \frac1{\sqrt n}\sum_{j=1}^n V_j(X_j)-\left(a+\frac{\b}{2\g}\right)\sqrt n\rb}\rb e^{-\frac{1}{4\g}s^2} ds}. 
		\end{align} 
		This completes the proof for the soft constraint $\eqref{2.31}$.
		
		To prove $\eqref{2.32}$ the idea is to use Proposition $\ref{prop:convergence}$, taking the limit $\g \to \infty$ of $\eqref{2.31}$. By convergence in distribution the characteristic function of $P_{n,\SC}$ converges pointwise to the characteristic function of $P_n$. We must first justify interchanging the limit and integration in the numerator and denominator on the r.h.s.~of $\eqref{2.31}$. It suffices to justify this only for the numerator as the denominator can be obtained by setting $t =0$. Defining $Z^j_{n,\b}:=\int_\R f_j(x_j) e^{-\b V_j(x_j)} dx$ and using independence under $\P_{n,\b}$, we get the bound
		\begin{align} 
			&\left|\E_{{n,\b}}\lb e^{is\lb \frac1{\sqrt n}\sum_{j=1}^n V_j(X_j)-\left(a+\frac{\b}{2\g}\right)\sqrt n\rb + it\langle \l,X-\E_{n,\b} X \rangle}\rb e^{-\frac{1}{4\g}s^2} \right| 
			\\ 
			&\leq\lvv \prod_{j=1}^{n} \int_{\R} \frac{1}{Z_{n,\b}^j} f_j(x_j) e^{-\b V_j(x_j)}  e^{is\frac{1}{\sqrt{n}}V_j(x_j)} e^{it\l_j x_j} dx_j\rvv\\
			& \leq C \prod_{j=1}^{n} \lvv  \int_{\R} f_j(x_j) e^{-\b V_j(x_j)}  e^{is\frac{1}{\sqrt{n}}V_j(x_j)} e^{it\l_j x_j} dx_j\rvv, \label{CFSPM}
		\end{align} 
		with $C:= \prod_{j=1}^n (Z^j_{n,\b})^{-1}< \infty$. The last inequality uses Assumption $\ref{ass:curve}$ which provides constants $\e,\d,M >0$ and an interval $J_j$ and such that
         \begin{align}
          {Z^j_{n,\b}} = \int_\R f_j(x_j) e^{-\b V_j(x_j)} dx_j \geq \int_{J_j} f_j(x_j) e^{-\b V_j(x_j)}dx_j > \e\d e^{-\b M } \text{ for all } j.\label{ZB}
          \end{align}
           We can now apply Lemma \ref{lem_vdc} to each of the integrals in the product in \eqref{CFSPM} for which $j \in G$, for $s$ sufficiently large, where we set $\t = t\l_j$, $r = s/\sqrt n$ and absorb the $\sqrt n$ term into the constant $C$. This poses no issue as $n$ is fixed. We can bound the factors for $j \not \in G$ by $K$ using Assumption \ref{ass:int} and write
          \begin{align}
              \left|\E_{{n,\b}}\lb e^{is\lb \frac1{\sqrt n}\sum_{j=1}^n V_j(X_j)-\left(a+\frac{\b}{2\g}\right)\sqrt n\rb + it\langle \l,X-\E_{n,\b} X \rangle}\rb \right| \leq C_1 |s|^{-\frac{|G|}{p}},
          \end{align}
          for some $C_1>0$ independent of $\g$. This bound paired with the trivial bound of 1 for $s$ small, is uniform in $\g$ and integrable in $s$ since $|G|>p$, so by Lebesgue's dominated convergence theorem we can interchange the limit $\g\to\infty$ and integration in the $s$-integral, with this limit now being trivial to evaluate. As mentioned above, the denominator is obtained setting $t=0$, and we note that this can be seen to be nonzero following the proof of Proposition \ref{prop:convergence}. This completes the proof of the proposition.
\end{proof}
\vspace{1em}
Finally we give the proof of Proposition \ref{lemma_beta} which ensures uniqueness and existence of the required $\b$.
\begin{proof}[Proof of Proposition \ref{lemma_beta}]
		We first deal with the hard constraint case. The constant $\b_0 \geq 0$ from \eqref{Mbound} of Assumption \ref{ass:BM} satisfies
        \begin{align}
            \E_{n,\b_0}\lb \sum_{j=1}^n V_j(X_j) \rb \geq na.
        \end{align}
        The idea is to show that for $\b$ large enough the opposite inequality holds, and then employing monotonicity and continuity to show the existence of a unique $\b$ satisfying the requirements. We first define $c : = a - \e /2$, where $\e$ denotes the constant of Assumption \ref{ass:mass}, and write
        \begin{align}
            &\E_{n,\b} (V_j(X_j)) = \frac{1}{Z^j_{n,\b}}\int_{V_j(x_j) \leq c} f_j(x_j) V_j(x_j) e^{-\b V_j(x_j)}  dx\\
            &+ \frac{1}{Z^j_{n,\b}}\int_{V_j(x_j) > c} f_j(x_j)V_j(x_j) e^{-\b V_j(x_j)} dx \leq c + \frac{1}{Z^j_{n,\b}}\int_{V_j(x_j) > c} V_j(x_j) f_j(x_j) e^{-\b V_j(x_j)}  dx.
        \end{align}
        On $\{V_j(x_j) > c\}$ we have $e^{-(\b -\b_0)V_j(x)} \leq e^{-(\b-\b_0)c}$. Using this and Assumptions \ref{ass:int} and \ref{ass:BM}, there exists $K>0$ such that
        \begin{align}
            &\int_{V_j(x_j) > c} f_j(x_j)V_j(x_j) e^{-\b V_j(x_j)} dx \leq  e^{-(\b-\b_0)c} \int_{V_j(x_j) > c} f_j(x_j)V_j(x_j) e^{-\b_0 V_j(x_j)} dx \\
            &\leq K e^{-(\b-\b_0)c},
        \end{align}
        for all $j$. From Assumption \ref{ass:mass} there exist constants $\e>0$, $\d>0$ and an interval $I_j$ with length $ |I_j| > \d$ and such that
        \begin{align}
            f_j(x_j) \geq \e \text{ and } e^{-\b V_j(x_j)} \geq e^{-\b(a-\e)}  \text{ for all }x \in I_j.
        \end{align}
        Therefore
        \begin{align}
            Z_{n,\b}^j \geq \int_{I_j} f_j(x_j) e^{-\b V_j(x_j)} dx \geq\e \d e^{-\b(a-\e)} \text{ for all }j.
        \end{align}
        Note that by \eqref{Mbound} of Assumption \ref{ass:BM}, as mentioned earlier, and H\"older's inequality,
\begin{align}
    a \leq \frac1n\sum_{j=1}^n\E_{n,\b_0}V_j(X_j)
    \leq \max_j\lb\E_{n,\b_0}V_j(X_j)^3\rb^{1/3}\leq K^{1/3},
\end{align}
so that $e^{\b_0 c}\leq e^{\b_0K^{1/3}}$.
        Combining this, there exists some $K_1>0$ such that
        \begin{align}
            \E_{n,\b} (V_j(X_j)) \leq a - {\e}/{2} + \frac{Ke^{\b_0 c}}{\e\d} e^{ -\b\e/2} = a - \e/2 +K_1e^{-\b \frac{\e}{2}},
        \end{align}
        for all $\b \geq \b_0$ and $j$. Since these constants are all independent of $j$ we can conclude that there exists a constant $\b^*$ uniform in $j$ such that for all $\b \geq \b^*$ we have
        \begin{align}
            \E_{n,\b}(V_j(X_j)) \leq a ,\text{ for all } j
        \end{align}
        and consequently
        \begin{align}
        \E_{n,\b}\lb \sum_{j=1}^n V_j(X_j)\rb \leq na, \label{bigbeta}
		\end{align} 
        holds for all $n$.
        To look at uniqueness, we observe by simple calculation that for all finite $\b\geq \b_0$ we have
		\begin{align}
			\frac{d}{d\b} \E_{{n,\b}}\left(\sum_{j=1}^n V_j(X_j)\right) = -\Var_{{n,\b}}\left(\sum_{j=1}^n V_j(X_j)\right) < 0, \label{vr}
		\end{align}
		where differentiation under the integral sign is justified by the integrability provided by Assumption \ref{ass:BM}. The inequality is strict since Assumption \ref{ass:curve} gives us $V_j$ non constant on the $J_j$, where $X_j$ has density bounded away from zero. In particular, $\b \mapsto \E_{n,\b} ( \sum_{j=1}^n V_j(X_j)) $ is continuous and strictly decreasing on $[\b_0,\infty)$. By \eqref{Mbound} of Assumption \ref{ass:BM} this map takes a value at least $na$ at $\b=\b_0$, and by $\eqref{bigbeta}$ it takes a value of at most $na$ at $\b = \b^*$. Hence there is a unique $\b \geq \b_0$ such that $\E_{{n,\b}}(\sum_{j=1}^n V_j(X_j)) = na $, and it lies in $ [\b_0,\b^*]$, proving the result for the hard constraint case.

        For the soft constraint case the argument is analogous, with the constant $na$ replaced by $n(a+\b/(2\g))$, which is strictly increasing in $\b$. The difference between the two sides of the required identity is therefore continuous and strictly decreasing on $[\b_0,\infty)$. It is nonnegative at $\b_0$ by \eqref{Mbound2} of Assumption \ref{ass:BM} and non-positive at $\b^*$ by \eqref{bigbeta} since
        \begin{align}
        	 \E_{n,\b^*}\lb \sum_{j=1}^n V_j(X_j)\rb \leq na \leq n \lb a + \frac{\b^*}{2\g}\rb.
        \end{align} 
        Hence there exists a unique $\b(\g) \in [\b_0,\b^*]$, with $\b^*$ independent of $\g$ and $n$.
	\end{proof}

\section{Auxiliary lemmas }
        This section prepares the ground for the proof of our main results Theorems $\ref{T1}$ and $\ref{T2}$.  We start by showing that the moment bounds given in Assumption \ref{ass:BM} for some $\b_0$, extend to any $\b \in [\b_0,\b^*]$, and thus to the $\b$ under which the constraint is fulfilled in expectation. 

 \begin{lemma} \label{lemma_b0}
    Let $\b^*$ be as in Proposition \ref{lemma_beta}. There exists a finite constant $\hat K$ such that
    for all $\b \in [\b_0,\b^*]$ we have
        \begin{align}
        \E_{n,\b} \lb V_j(X_j)^3 \rb \leq \hat K, \quad \E_{n,\b} \lb |X_j|^3 \rb \leq \hat K, \quad \Var_{n,\b}(V_j(X_j)) \leq \hat K, \quad \Var_{n,\b}(X_j) \leq \hat K
        \end{align}
for all $j$ and $n$.
\end{lemma}
\begin{proof}
    We prove that these bounds hold for $V_j(X_j)$, showing the same for $X_j$ is analogous. By Assumptions \ref{ass:int} and \ref{ass:BM}  there exists $K>0$ such that
    \begin{align}
        \E_{n,\b} \lb V_j(X_j)^3 \rb = \frac{Z^j_{n,\b_0}}{Z^j_{n,\b}} \E_{n,\b_0} \lb V_j(X_j)^3 e^{-(\b -\b_0) V_j(X_j)} \rb 
        \leq \frac{Z^j_{n,\b_0}}{Z^j_{n,\b}} \E_{n,\b_0} \lb V_j(X_j)^3 \rb \leq  {Z^j_{n,\b}}^{-1} K^2 \label{star1}
    \end{align}
    holds for all $j$, where we have used $\b \geq \b_0$ and $V_j \geq 0$.  Using that $\b \leq \b^*$ we can conclude from $\eqref{ZB}$ there exists $\e,\d>0$ such that
    \begin{align}
        {Z^j_{n,\b}} >\e\d e^{-\b M}\geq \e\d e^{-\b^* M}  \text{ for all }j. \label{star}
    \end{align}
    Now since all constants in \eqref{star1} and \eqref{star} are uniform in $j$ combining these bounds gives the required bound for the expectation.
    \newline \indent
        The variance bound follows immediately by H\"older's inequality which gives
    \begin{align}
    \Var_{n,\b}(Y) \leq \E_{n,\b} \lb Y^2 \rb \leq \lb\E_{n,\b} \lb |Y|^3\rb \rb^{\frac{2}{3}},
    \end{align}
   applied with $Y = V_j(X_j)$, again using $V_j\geq0$.
    
\end{proof}
As our representation result Proposition \ref{prop_CF} shows, we have to understand the interplay between the two statistics $\sum_{j=1}^n\l_j(X_j-\E_{n,\b}X_j)$ and $\frac1{\sqrt n}\sum_{j=1}^n V_j(X_j)-a\sqrt n$. In order to relate the different scalings of the two statistics, namely $\l_j$ and $1/\sqrt n$, it will be convenient to define $\hat \l : = \l \sqrt n$. The first statistic is linear while the second is not linear everywhere thanks to the curvature of $S_{a,V}$ not being 0 everywhere by Assumption \ref{ass:curve}. The next lemma shows that for each $j$ the variance of any combination $sV_j(X_j)+t\hat\l_jX_j$ is bounded below by a multiple of $s^2 +t^2\hat \l_j^2$, uniformly in $j$. This will translate into a bound on the decay of the characteristic function of Proposition \ref{prop_CF} in Proposition \ref{lemma_CF_unif} below.

\begin{lemma} \label{lemma_mp}
Defining $\hat \l := \l \sqrt n$, there exist constants $0 <  \e \leq C < \infty$ such that for all $\b\in [\b_0,\b^*]$, $t \in \R$, $s\in\R$,
	\begin{align}
		\e(s^2+t^2||\l||_2^2) \leq \frac{1}n \sum_{j=1}^n \Var_{n,\b} \lb sV_{j}(X_j) + t\hat\l_jX_j\rb \leq C (s^2 +t^2||\l||_2^2).
	\end{align}
\end{lemma}
\begin{proof}
	Fix $s,t \in \R$ and for each $j$ set
	\begin{align}
		g_j(x) : = sV_j(x) + t \hat \l_j x, \quad \rho_j := \sqrt{ s^2 + t^2 \hat \l_j^2 } \geq 0.
	\end{align}
	We will prove the coordinate-wise bound 
	\begin{align}
		\e \rho_j^2 \leq \Var_{n,\b}(g_j(X_j)) \leq C \rho_j^2, \label{ineq}
	\end{align}
	with $\e,C$ independent of $j,s$ and $t$, and the lemma follows by averaging since
	\begin{align}
		\frac{1}{n} \sum_{j=1}^n \rho_j^2 = s^2 + t^2\lb \frac{1}{n} \sum_{j=1}^n \hat \l_j ^2 \rb = s^2 + t^2 ||\l||^2_2.
	\end{align}
	For the upper bound we can write
	\begin{align}
		\Var_{n,\b} (g_j(X_j)) \leq 2 s^2 \Var_{n,\b} (V_j(X_j))  + 2t^2 \hat \l_j ^2 \Var_{n,\b}(X_j) \leq C\rho_j^2
	\end{align}
	for some $C>0$ uniform in $j$ and $n$, where we can use Lemma \ref{lemma_b0} to conclude finiteness of $C$. \\
		\indent 
        For the lower bound we may assume $\rho_j \neq 0$, as this case is trivial. By Assumption $\ref{ass:curve}$ there exist constants $\e_0>0$, $\d_0>0$, independent of $j$, and an interval $J_j \subset \R \setminus H_j$ with $|J_j| > \d_0$ such that 
	\begin{align}
		|V''_j(x)| > \e_0,\quad V_j(x) \leq M, \quad f_j(x) > \e_0,\quad \text{ for } x \in J_j. \label{nice}
	\end{align}
    For simplicity we set $\d:= \d_0/4$, so that $|J_j| > 4\d$. Since $V_j''$ is continuous and bounded away from zero on $J_j$ it has constant sign on $J_j$. Hence $V_j$ is convex or concave on $J_j$. Let $(a_j,b_j)$ denote the interval $J_j$. If $V_j$ is convex on $J_j$, the $V'_j$ is non-decreasing, so for $x \in J_j$ we can write
    \begin{align}
        -\frac{M}{x-a_j} \leq \frac{ V_j(x) - V_j(a_j)}{x-a_j} \leq V'_j(x) \leq \frac{ V_j(b_j) - V_j(x)}{b_j-x} \leq \frac{M}{b_j-x}.
    \end{align}
    If $V_j$ is concave on $J_j$, then $V'_j$ is non-increasing, so analogously we have for $x \in J_j$
    \begin{align}
        -\frac{M}{b_j-x} \leq  V'_j(x) \leq  \frac{M}{x-a_j}.
    \end{align}
    In either case $|V'_j(x)| \leq M/ (\min\{x-a_j,b_j-x\})$. Letting $\tilde J_j := (a_j+\d,b_j-\d)$ we have $|\tilde J_j| > 2\d$ and every point of $\tilde J_j$ is at least a distance $\d$ away from both $a_j$ and $b_j$. Hence
    \begin{align}
        |V''_j(x)| > \e_0, \quad V_j(x) \leq M, \quad |V'_j(x)| \leq L : = \frac{4M}{\d},\quad  f_j(x) > \e_0 \text{ for all } x \in \tilde J_j.
    \end{align}
    Since $0 \leq \b \leq \b^*$ and $V_j \geq 0$ it follows that
    \begin{align}
        \tilde f_j(x) : = f_j(x) e^{-\b V_j(x)} \geq \e_0 e^{-\b^* M} =: \hat \e>0 \text{ on } \tilde J_j.
    \end{align}
	Defining $\hat C := (2 \sup_j (Z^j_{n,\b})^2)^{-1}$, Assumption \ref{ass:int} and \eqref{star} show that $\hat C$ is non-zero and finite. Recall that if $v$ is a probability measure on $\R$, $X \sim v$, and $h : \R \to \R$ is measurable with $\E(h(X)^2) < \infty$, then
	\begin{align}
		\Var(h(X)) = \frac{1}{2} \int \int (h(x)-h(y))^2 dv(x)dv(y).
	\end{align}
	We apply this with $v$ being the distribution of $X_j$ under $P_{n,\b}$ and $h = g_j$ which is square integrable by Lemma \ref{lemma_b0}. Since the integrand is nonnegative we may restrict to $\tilde J_j \times \tilde J_j$, where $\tilde f_j \geq \hat \e$ to write,
	\begin{align}
		\Var_{n,\b}(g_j(X_j)) \geq \hat C \hat\e^2 \int_{\tilde J_j} \int_{\tilde J_j} (g_j(x) - g_j(y) )^2 dx dy.
	\end{align}
	If $|g'_j| \geq \rho>0$ on some interval $I \subset \tilde J_j$, then $g'_j$ being continuous and nonvanishing, has constant sign on $I$ and the  mean value theorem gives us $|g_j(x)-g_j(y)| \geq \rho |x-y|$ for $x,y \in I$. Hence
	\begin{align}
		\Var_{n,\b}(g_j(X_j)) \geq \frac{\hat C \hat \e ^2\rho^2|I|^4}{6}. \label{killer}
	\end{align}
	We now fix 
	\begin{align}
		\theta : = \min \lb \frac{1}{\sqrt 2}, \frac{1}{2\sqrt 2 L}\rb > 0,
	\end{align}
	and distinguish two cases according to the size of $|s|$ relative to $\rho_j$.
	\\
	\indent 
	The first case is $|s| \leq \theta\rho_j$. Then $t^2\hat \l_j^2 = \rho_j^2 -s^2 \geq (1-\theta^2)\rho_j^2 \geq \rho_j^2 /2$, so $|t\hat \l_j| \geq \rho_j/\sqrt 2$. For $x \in \tilde J_j$,
	\begin{align}
		|g'_j(x)| \geq |t||\hat\l_j| - |s||V_j'(x)| \geq \frac{\rho_j}{\sqrt 2} - \theta L \rho_j \geq  \frac{\rho_j}{2\sqrt 2} = \rho_j \rho_1, \quad \rho_1 := \frac{1}{2\sqrt 2},
	\end{align}
	and $\eqref{killer}$ with $I = \tilde J_j$ gives
	\begin{align}
		\Var_{n,\b}(g_j(X_j)) \geq \frac{\hat C \hat \e^2 \rho_j^2 \rho_1^2 (2\d)^4 }{6}. \label{BND1}
	\end{align}
	\indent
	In the second case we have $|s| \geq \theta\rho_j$. Here $|g''_j(x)| = |s| |V''_j(x)| \geq \theta \e_0\rho_j$. Since $V''_j$ is continuous on $\tilde J_j$ and nowhere zero, it has constant sign. Therefore $g'_j$ is strictly monotonic and has at most one zero $\hat x_j$ on $\tilde J_j$. A Taylor expansion around $\hat x_j$ gives 
	\begin{align}
		|g'_j(x)| \geq \theta \rho_j\e_0 |x-\hat x_j| \text{ on } \tilde J_j.
	\end{align}
	Since $|\tilde J_j| > 2\d$, the point $\hat x_j$ is at least $\d$ from one endpoint of $\tilde J_j$. If no such $\hat x_j$ exists, take $\hat x_j$ to be the endpoint of $\tilde J_j$ for which $|g'_j|$ acquires its minimum in which case $g'_j$ has constant sign and the same bound follows from monotonicity. Letting $ \hat I_j$ be the sub interval of length $\d/2$ on that side, at distance at least $\d/2$ away from $\hat x_j$, we have 
	\begin{align}
		|g_j'| \geq \frac{\theta \rho_j\e_0 \d}{2} = \rho_j\rho_2 \text { on } \hat I_j, \quad \rho_2 : = \frac{\theta \e_0 \d}{2}.
	\end{align}
	Applying $\eqref{killer}$ with $I = \hat I_j$ gives 
	\begin{align}
		\Var_{n,\b} (g_j(X_j)) \geq \frac{\hat C \hat \e^2 \rho_j^2\rho_2^2(\d/2)^4}{6} \label{BND2}
	\end{align}
	for all $j$.\\
	\indent
	In either case the bound holds for every $j$, and all constants are uniform in $j$, $s$ and $t$. Hence, \eqref{ineq} holds with 
	\begin{align}
		\e : = \frac{ \hat C \hat \e^2}{6}  \min \lb \rho_1^2(2\d)^4, \rho_2^2(\d/2)^4\rb >0,
	\end{align}
	and averaging over $j$ completes the proof.
\end{proof}

\vspace{1em}
The next lemma shows that the variance terms in Theorems \ref{T1} and \ref{T2} are bounded away from $0$ and infinity. 
\medskip

\begin{lemma} \label{sigmalemma}
	Let $\s_n$ be defined as in Theorem \ref{T1} or Theorem \ref{T2}. Then under the conditions of these theorems we have
	\begin{align}
		0< \liminf_{n \to \infty} \s_n^2 \leq \limsup_{n \to \infty} \sigma_n^2 < \infty.
	\end{align}
\end{lemma}
\begin{proof}
	We first consider the hard constraint case (Theorem \ref{T1}), and the soft constraint case is treated at the end of the proof. Let
	\begin{align}
		Y := \lb \sum_{j=1}^n V_{j}(X_j), \sum_{j=1}^n \l_j X_j \rb, \quad M_n := \Cov_{n,\b}(Y) = \begin{pmatrix}
			A&q\\q&d
		\end{pmatrix},
	\end{align}
	with 
	\begin{align}
		A:=  \sum_{j=1}^n \Var_{n,\b} (V_{j}(X_j)),\quad q := \sum_{j=1}^n  \l_j \Cov_{{n,\b}}(X_j,V_j(X_j)) ,\quad  d := \sum_{j=1}^n {\l_j}^2\Var_{{n,\b}}(X_j). 
	\end{align}
	Then $\s_n^2 = d -q^2/A$ is the Schur complement of $A$ in $M_n$. 
	For $v = (v_1,v_2) \in \R^2$, the independence structure of $P_{n,\b}$ gives
	\begin{align}
		v^T M_n v = \sum_{j=1}^n \Var_{n,\b}(v_1 V_j(X_j) + v_2 \l_j X_j).
	\end{align}
	Writing $\l_j = \hat\l_j/\sqrt n$ and applying Lemma $\ref{lemma_mp}$, which applies since $\b \in [\b_0,\b^*]$ by Proposition \ref{lemma_beta}, with $s=v_1$ and $t = v_2/\sqrt n$, yields
	\begin{align}
		\e(nv_1^2 + v_2^2||\l||^2_2) \leq v^TM_n v \leq C(nv_1^2+v^2_2||\l||^2_2). \label{circle}
	\end{align}
	Taking $v = (1,0)$ in the lower bound gives $A \geq \e n >0$, so $A$ is invertible and $\s^2_n$ is well defined. Completing the square, for every $w \in \R$ we have
	\begin{align}
		\begin{pmatrix}
			w\\1
		\end{pmatrix}^T M_n \begin{pmatrix}
			w\\1
		\end{pmatrix} = Aw^2+2qw+d = A\lb w+\frac{q}{A}\rb^2 + \s^2_n,
	\end{align}
	hence
	\begin{align}
		\s^2_n = \min_{w \in \R} \begin{pmatrix}
			w\\1
		\end{pmatrix}^T M_n \begin{pmatrix}
			w\\1
		\end{pmatrix}.
	\end{align}
	Applying $\eqref{circle}$ with $v = (w,1)$ gives 
	\begin{align}
		\e(nw^2+||\l||^2_2) \leq \begin{pmatrix}
			w\\1
		\end{pmatrix}^T M_n \begin{pmatrix}
			w\\1
		\end{pmatrix} \leq C(nw^2+||\l||^2_2),
	\end{align}
	and minimizing over $w$ gives 
	\begin{align}
		\e||\l||^2_2 \leq \s_n^2 \leq C||\l||^2_2.
	\end{align}
    \indent
    In the soft constraint case, $A$ is replaced by $A + n/(2\g) \geq A$. Hence $d-q^2/(A+n/(2\g))$
    lies between $d-q^2/A$ and $d$, and the same bounds hold. This concludes the proof using the assumptions on $\l$ as in \eqref{LASS}, which apply both for Theorems \ref{T1} and \ref{T2}. 
\end{proof}

\vspace{1em}
The final proposition of this section provides a uniform asymptotic bound for the joint characteristic function of $V_j(X_j)$ and $X_j$ under $P_{n,\b}$. This bound will be important when employing Proposition \ref{prop_CF} to prove Theorem \ref{T1}. More precisely, the bound will allow us to interchange the limit $n \to \infty$ and the $s$ integral in representation \eqref{2.32} of Proposition \ref{prop_CF} (see \eqref{CF} below).
\medskip

\begin{prop} \label{lemma_CF_unif}
	For every $K>0$ there exist constants $c>0$, $C>0$ and $C'< \infty$, uniform in $j$ and $n$ such that the following holds. With $\b$ as in Proposition $\ref{lemma_beta}$, let $q \in \R$ and let $\hat t \in \R^n$ satisfy  $||\hat t||_2 \leq K$. Then writing  
	\begin{align}
		\phi_j(q) := \E_{n,\b} \left( e^{iqV_j(X_j) + i \hat t_j X_j}\right),
	\end{align}
	we have
	\begin{align}
		\prod_{j=1}^n\lvv\phi_j(q)\rvv \leq C' (1+c\min(|q|,|q|^2))^{-Cn}.
	\end{align}
\end{prop} 
\begin{proof} 
	For the convenience of the reader we start by giving a broad outline of the rather lengthy proof. We will distinguish three major cases according to the size of $\lv q\rv$. For $\lv q\rv$ small enough, a Taylor expansion of $\phi_j(q)$ around $q=0$ will provide the sought bound. For $\lv q\rv$ large enough, Lemma \ref{lem_vdc} will give the desired. For $\lv q\rv$ in the intermediate range it will suffice to show $\lv \phi_j(q)\rv\leq c<1$ for some $c>0$ uniform in $j$. This bound will be shown using that $|\phi| = 1$ forces a lattice distribution. The details of the proof take up considerable space as these outlined ideas can not be employed directly, in particular due to $\phi_j$ being a joint rather than a single characteristic function.
	
	We now begin the actual proof with the case $\lv q\rv<\e_0$, for some small $\e_0>0$ independent of $j$. An expansion around $q=0$ then yields
	\begin{align}
		\phi_j\lb q \rb=
		\E_{n,\b}\left(e^{i\hat t_j X_j}\right)  + qi\E_{n,\b}\left(V_j(X_j)e^{i\hat t_j X_j}\right) -\frac{q^2}{2}\E_{n,\b}\left(V_j(X_j)^2e^{i\hat t_j X_j}\right) + \O\lb |q|^3 \rb,  \label{q-expansion}
	\end{align}
	where the $\O$-term is uniform in $j$ by the moment bounds of Lemma $\ref{lemma_b0}$. Using \eqref{q-expansion}, we get
	\begin{align}
		&\lvv\phi_j\lb q \rb\rvv^2 =\phi_j\lb q \rb \overline{\phi_j\lb q \rb}\\
		&=\E_{n,\b}\left(e^{-i\hat t_j X_j}\right)  \E_{n,\b}\left(e^{i\hat  t_j X_j}\right)\\
		&+qi\left[ \E_{n,\b}\left(e^{-i\hat t_j X_j}\right)\E_{n,\b}\left(V_j(X_j)e^{i\hat t_j X_j}\right) -\E_{n,\b}\left(e^{i\hat t_j X_j}\right)\E_{n,\b}\left(V_j(X_j)e^{-i\hat t_j X_j}\right)\right]\\
		&+q^2\E_{n,\b}\left(V_j(X_j)e^{-i\hat t_j X_j}\right)\E_{n,\b}\left(V_j(X_j)e^{i\hat t_j X_j}\right)\\
		&-\frac{q^2}2\left[\E_{n,\b}\left(e^{i\hat t_j X_j}\right)\E_{n,\b}\left(V_j(X_j)^2e^{-i\hat t_j X_j}\right) +\E_{n,\b}\left(e^{-i\hat t_j X_j}\right)\E_{n,\b}\left(V_j(X_j)^2e^{i\hat t_j X_j}\right)\right]+\O\left(|q|^3\right),
	\end{align} 
	with the order term being uniform in $j$.
	For the order $q^0$ term we use the crude bound
	\begin{align}
		\lvv \E_{n,\b}\left(e^{-i\hat t_j X_j}\right)  \E_{n,\b}\left(e^{i\hat  t_j X_j}\right)\rvv \leq 1.
	\end{align}
	Now expanding in $\hat t_j$ around 0, we get 
	\begin{align}
		&\E_{n,\b}\left(V_j(X_j)e^{-i\hat t_j X_j}\right)\E_{n,\b}\left(V_j(X_j)e^{i\hat t_j X_j}\right)-\frac{1}2\left[\E_{n,\b}\left(e^{i\hat t_j X_j}\right)\E_{n,\b}\left(V_j(X_j)^2e^{-i\hat t_j X_j}\right)\right] \\
		&-\frac{1}2 \left[\E_{n,\b}\left(e^{-i\hat t_j X_j}\right)\E_{n,\b}\left(V_j(X_j)^2e^{i\hat t_j X_j}\right)\right] =  - \Var_{n,\b}(V_j(X_j)) + \O \lb \hat t_j \rb.
	\end{align}
	for the order $q^2$ term. 
	Performing this expansion for the remaining term we write
	\begin{align}
		&\lvv \phi_j\lb q \rb\rvv^2 \leq 1 - q^2\Var_{n,\b}(V_j(X_j)) + \O(|q||\hat t_j|)+\O(|q|^3),
	\end{align}
	with all order terms giving constants uniform in $j$. Here the contribution $q^2\O(|\hat t_j|)$ from the expansion of the $q^2$ coefficients is included in $\O(|q||\hat t_j|)$, since $|q| \leq \e_0 \leq 1$. By the coordinate wise bound \eqref{ineq} in the proof of Lemma \ref{lemma_mp}, applied with $s = 1,t=0$, there exists $\e_1>0$ such that 
	\begin{align}
		\Var_{n,\b}(V_j(X_j)) > \e_1,
	\end{align}
	for all $j$. We can absorb the $\O(|q|^3)$ into the $|q|^2$ term by making $\e_0$ sufficiently small and write
	\begin{align}
		|\phi_j(q)|^2 \leq 1 - |q|^2 \frac{\e_1}{2} + T_0 |q| |\hat t_j|,
	\end{align}
	for some constant $T_0>0 $ uniform in $j$. Using Young's inequality we can derive
	\begin{align}
		T_0 |q| |\hat t_j| \leq \frac{\e_1}{4}|q|^2 + \frac{T_0^2}{\e_1} \hat t_j^2.
	\end{align}
	Hence, there exist constants $\hat \e>0$ and $T>0$ such that 
	\begin{align}
		|\phi_j(q)|^2 \leq 1 - |q|^2 \hat \e + T \hat t_j^2
	\end{align}
	for all $|q| \leq \e_0$. Using $1+x \leq e^x$, we obtain
	\begin{align}
		\prod_{j=1}^n |\phi_j(q)|  \leq  \exp \lb \frac12 \sum_{j=1}^n (-|q|^2\hat\e + T\hat t_j^2)\rb.
	\end{align}
	Using $||\hat t||_2^2 \leq K^2$, we have
	\begin{align}
		\prod_{j=1}^n |\phi_j(q)| \ \leq \exp\lb\frac{TK^2}{2}\rb \exp \lb -|q|^2\hat\e \frac{n}2\rb.
	\end{align}
	Hence there exists $C'=C'(K)>0$ and $\hat C_1>0$, $c_1>0$ independent of $K$ such that 
	\begin{align}
		\prod_{j=1}^n |\phi_j(q)| \leq C'(1+c_1 |q|^2)^{-\hat C_1 n} \label{c1}
	\end{align}
	for all $|q| \leq \e_0$.\\
	\indent  We now move to the case of large $q$. Here we can apply Lemma \ref{lem_vdc} with $r = q$, $\t = \hat t_j$ to all $j \in G$,  which gives 
    \begin{align}
         \lvv \phi_j\lb q \rb\rvv =  \lvv \int_\R \frac{1}{Z^j_{n,\b}} e^{-\b V_j(x)+i\hat t_j x} f_j(x) e^{iqV_j(x)}dx \rvv \leq C C_0(1+\b + |\hat t_j|)|q|^{-1/p} \leq C_1|q|^{-1/p}, \label{boundq}
    \end{align}
    for some $C_1>0$ independent of $j$ since $\b \leq \b^*$ and $\|\hat t\|_{\max} \leq K$, where we have bounded $Z^j_{n,\b}$ using \eqref{ZB}. Using the trivial bound $|\phi_j(q)| \leq 1$ for all $j \not \in G$, and the bound \eqref{boundq} for $j \in G$ with $|q| \geq k_0$ for sufficiently large $k_0$ we can write
    \begin{align}
        \prod_{j=1}^n |\phi_j(q)| \leq (1+c_2|q|)^{-\hat C_2 n}, \label{c2}
    \end{align}
    for some $c_2>0$ and $\hat C_2 >0$, utilizing $|G| > cn$ for some $c>0$ independent of $n$.
     \\
	\indent
	We now move to the case $|q| \in [\e_0,k_0]$ for some $\e_0,k_0>0$. We first make the following change of notation, defining 
	\begin{align}
		Y_{j,q} := qV_j(X_j) + \hat t_j X_j.
	\end{align}
	We then have
	\begin{align}
		\phi_j(q) = \E_{n,\b} \left( e^{iqV_j(X_j) + i \hat t_j X_j}\right) = \E_{n,\b}\left(e^{iY_{j,q}} \right) = \phi_{Y_{j,q}}(1),
	\end{align}
	where $\phi_{Y_{j,q}}$ is the characteristic function of $Y_{j,q}$ under $P_{n,\b}$ evaluated at one.
	Since
	\begin{align}
		 \E_{n,\b}\left(\lvv V_j(X_j)\rvv\right) \leq \hat K^{\frac{1}{3}} \text{ and } \E_{n,\b}\left(\lvv X_j\rvv\right) \leq  \hat K^{\frac{1}{3}} , 
	\end{align}
	by Lemma $\ref{lemma_b0}$, $ \lv q\rv \in [\e_0,k_0]$ and $||\hat t||_2 \leq K$ by assumption, the family of random variables $Y_{j,q}$ indexed by $n,j,q \in[\e_0,k_0]$ and $\hat t$ with $||\hat t ||_2 \leq K$ is tight. With this in mind, we claim that there is $\bar c<1$ such that 
    \begin{align}
         |\phi_j(q)| \leq \bar c \label{NM}
    \end{align}
   for all $n$, $j$, $|q| \in[\e_0,k_0]$.
	If this was not true, there must exist $n_m,j_m,q_m$ and $\hat t^{(m)}$ in this range with $|\phi_{Y_m}(1)| \to 1$, where $Y_m := Y_{j_m,q_m}$. By tightness and Helly's Selection Theorem, a subsequence $(Y_m)$ then converges in distribution to some random variable $Y$. By equicontinuity of characteristic functions of a tight sequence of random variables if such a sequence $m$ exists we must have
	\begin{align}
		\lvv \phi_{Y}(1) \rvv = \lim_{m \to \infty} \lvv \phi_{Y_m}(1)\rvv = 1,
	\end{align}
    where the limit is taken along this subsequence. It is however well-known that for a characteristic function $\phi$ the property $\lv \phi(t)\rv=1$ for some real $t\not=0$ implies that the underlying distribution is a lattice distribution (cf. \cite[Theorem 2.1.4]{Lukacs}). Hence $Y$ has a lattice distribution.
	We now examine the density of the random variable 
	\begin{align}
		Y_{j,q} = qV_j(X_j) + \hat t_j X_j,
	\end{align}
	to show that this is not possible, thus proving the claim $\eqref{NM}$.
	We again define
	\begin{align}
		g_j(x_j) := qV_j(x_j) + \hat t_j x_j
	\end{align} 
	and analogous to the second case of the proof of Lemma \ref{lemma_mp}, using $|g''_j|$ bounded below on $J_j$ the interval supplied by Assumption \ref{ass:curve}, there exists a sub interval $\tilde I_j$ with $|\tilde I_j| > \d$ for some $\d >0$ uniform in $j$ and we have 
	\begin{align}
		|V_j'(x_j)| \leq L, \quad e^{-\b V_j(x_j)}f_j(x_j) \geq \hat\e  \text{ and } |g'_j(x_j)| \geq \hat \e \text{ for } x \in \tilde I_j, \label{keybound}
	\end{align}
	with $\hat \e >0$ and $L>0$ uniform in $j$. Using this we can write $p_{Y_{j,q}}(y_j)$, the density of the law of $Y_{j,q}$ restricted to $\{ X_j \in \tilde I_j\}$ as
	\begin{align}
		p_{Y_{j,q}}(y_j) = \frac1{Z_{n,\b}^j} \frac{1}{\lvv q V'_j(g_j^{-1}(y_j)) + \hat t_j\rvv} f_j(g_j^{-1}(y_j))e^{-\b V_j(g_j^{-1}(y_j))}, \text{ for } y_j \in  g_j(\tilde I_j).
	\end{align} 
	Using \eqref{keybound} and $||\hat t||_2 \leq K$, we can conclude that this density is bounded below by a constant $c_4>0$ on this interval $g_j(\tilde I_j)$, with the length of such interval $|g_j(\tilde I_j)|> \hat\e\d$, with $c_4,\hat\e$ and $\hat \d$ uniform in $n,j,q \in [\e_0,k_0]$ and $\hat t$. We can conclude that any subsequence limit of $Y_{j,q}$ cannot be supported on a lattice, and the claim $\eqref{NM}$ holds for $|q| \in [\e_0,k_0]$. We now have for $|q| \in [\e_0,k_0]$, for any $\e_0,k_0 >0$ there exists $\hat C_3>0$, $c_3>0$ such that
	\begin{align}
		\prod_{j=1}^n \lvv \phi_j(q) \rvv \leq (1+c_3|q|)^{-\hat C_3 n}. \label{c3} 
	\end{align}
	Now by taking $c := \min(c_1,c_2,c_3)$, $C := \min (\hat C_1,\hat C_2,\hat C_3)$, and $C'' = \max(C',1)$ from each of the bounds $\eqref{c1},\eqref{c2},\eqref{c3}$, we have
	\begin{align}
		\prod_{j=1}^n \lvv \phi_j(q) \rvv \leq C''(1+c\min(|q|,|q|^2))^{-{Cn}}
	\end{align}
	for all $q$, with $c$ being independent of $j$.
\end{proof}

\section{Proof of Main Results}
	We are now in a position to prove our two main results Theorem $\ref{T1}$ and Theorem $\ref{T2}$. We take $\b$ as determined by Proposition $\ref{lemma_beta}$, for the soft and hard constraints, respectively. We first prove Theorem $\ref{T2}$, the soft constraint case, the hard constraint case will be discussed afterwards.

	\bigskip 
	\begin{proof}[Proof of Theorem \ref{T2}]
		The proof relies on the key structure of the characteristic function of $P_{n,\SC}$, which has been developed in Proposition $\ref{prop_CF}$ given by
		\begin{align}
			&\E_{{n,\SC}}\lb e^{it\s_n^{-1}\langle \l,X-\E_{n,\b} X \rangle}\rb \\
			&= \frac{\int_{\R}\E_{{n,\b}}\lb e^{is\lb \frac1{\sqrt n}\sum_{j=1}^n V_j(X_j)-\left(a+\frac{\b}{2\g}\right)\sqrt n\rb + it\s_n^{-1}\langle \l,X-\E_{n,\b} X \rangle}\rb e^{-\frac{1}{4\g}s^2} ds}{\int_{\R}\E_{{n,\b}}\lb e^{is\lb \frac1{\sqrt n}\sum_{j=1}^n V_j(X_j)-\left(a+\frac{\b}{2\g}\right)\sqrt n\rb}\rb e^{-\frac{1}{4\g}s^2} ds}.\label{CF2}
		\end{align}
		We will show that this characteristic function converges as $n \to \infty$ to that of a normal random variable. To do this we wish to take the limit $n \to \infty$ inside the $s$ integrals of $\eqref{CF2}$, justifying this by use of Lebesgue's dominated convergence theorem. From there we can evaluate this limit for fixed $s,t$ with the Lindeberg Feller CLT, utilizing the independence structure of $P_{n,\b}$. Finally we evaluate these $s$ integrals and show that what remains is the characteristic function of a standard normal distribution. \newline
		\indent We first discuss the asymptotics of the numerator of $\eqref{CF2}$, the asymptotics of the denominator can be obtained by setting $\l = 0$. Examining the integral of the numerator of $\eqref{CF2}$ for $s$ fixed, and applying the trivial bound
		\begin{align}
			\lvv\E_{{n,\b}}\lb e^{is\lb \frac1{\sqrt n}\sum_{j=1}^n V_j(X_j)-\left(a+\frac{\b}{2\g}\right)\sqrt n\rb + it\s_n^{-1}\langle \l,X-\E_{n,\b} X \rangle}\rb e^{-\frac{1}{4\g}s^2}\rvv  \leq   e^{-\frac{s^2}{4\g}}   \label{Gausdecay}
		\end{align}
		allows, by use of Lebesgue's dominated convergence theorem, for the limit $n\to\infty$ to be brought inside the integral. Recalling $\hat \l = \l \sqrt n$, and defining
		\begin{align}
			Y_j := sV_j(X_j) + \frac{t\hat{\l}_j X_j}{\sigma_n} - \E_{n,\b}\lb sV_j(X_j) + \frac{t\hat{\l}_j X_j}{\sigma_n} \rb 
		\end{align}	
		we can use
        \begin{align}
            \E_{n,\b} \lb \sum_{j=1}^n V_j(X_j) \rb =n\lb a+\frac{\b}{2\g}\rb
        \end{align} to write
		\begin{align}
			\E_{{n,\b}}\lb e^{is\lb \frac1{\sqrt n}\sum_{j=1}^n V_j(X_j)-\left(a + \frac{\b}{2\g}\right)\sqrt n\rb + it\s_n^{-1}\langle \l,X-\E_{n,\b} X \rangle}\rb = \E_{{n,\b}} \left(e^{\frac{i}{\sqrt{n}}\sum_{j=1}^{n}Y_j} \right) \label{Yjsetup},
		\end{align}
		Under $P_{n,\b}$ the random variables $Y_j$ are centered, independent and  non-identically distributed. Having justified bringing the limit $n \to \infty$ inside the integral, we are in a suitable position to apply the Lindeberg Feller CLT. We first define $s_n^2 := \sum_{j=1}^{n}\Var_{{n,\b}} Y_j$. For $(s,t) \neq (0,0)$, Lemma \ref{lemma_mp}, with $t/\s_n$ in place of $t$, and Lemma \ref{sigmalemma} give 
		\begin{align}
			\d < \frac{s_n^2}{n} < \d^{-1} 
		\end{align}
		for all $n$ sufficiently large, with $\d>0$ depending only on $s$ and $t$. For $s=t=0$ there is nothing to prove.
        \newline \indent 
        Now the Lindeberg Feller CLT tells us 
		\begin{align}
			\E_{{n,\b}} \left(\exp\left(\frac{i}{\sqrt n}\sum_{j=1}^{n}Y_j\right)\right)  = 	\exp\left( -\frac{s_n^2}{2n} \right) + o(1) \label{Lind}
		\end{align}
		under the Lindeberg Feller condition, which reads
		\begin{align}
			\lim_{n \to \infty} \frac{1}{n} \sum_{j=1}^n \E_{n,\b}\lb Y_j^2 \dopp1_{\{|Y_j| >\e\sqrt n\}} \rb = 0 \text{ for all }\e > 0.
		\end{align}  
		We note that the exponential term on the right hand side of \eqref{Lind} is dominated by an integrable function from Lemma $\ref{lemma_mp}$ and hence poses no obstacle to our interchanging of limits.
		To see that this Lindeberg Feller condition holds, we begin by writing
		\begin{align}
			\frac{1}{n} \sum_{j=1}^n \E_{n,\b} \lb Y_j^2 \dopp1_{\{|Y_j| >\e\sqrt n\}} \rb \leq \frac{1}{\e n^{3/2}} \sum_{j=1}^n \E_{n,\b} \lb |Y_j|^3 \rb.
		\end{align}
		Then using convexity and Lemma \ref{lemma_b0}, there exists $K>0$ such that
		\begin{align}
			\sum_{j=1}^n \E_{n,\b}\lb |Y_j|^3 \rb \leq K \lb n |s|^3 + \frac{|t|^3}{\s_n^3} \sum_{j=1}^n |\hat \l_j^3 |\rb.
		\end{align}
		With $\hat\l_j = \sqrt n \l_j$ we have 
		\begin{align}
			n^{-3/2} \sum_{j=1}^n |\hat \l_j|^3 = \sum_{j=1}^n |\l_j|^3 = ||\l||^3_3,
		\end{align}
		and therefore,
		\begin{align}
			\frac{1}{n} \sum_{j=1}^n \E_{n,\b} \lb Y_j^2 \dopp1_{\{|Y_j| >\e\sqrt n\}} \rb \leq \frac{K}{\e} \lb \frac{|s|^3}{\sqrt n} + \frac{|t|^3}{\s_n^3} ||\l||^3_3 \rb.
		\end{align}
		The first term tends to zero trivially and to see that the second term goes to zero we can write $||\l||^3_3 \leq ||\l||_{\max} ||\l||^2_2$ noting that \eqref{LASS} ensures this tends to zero, and bound $\s_n$ from below using Lemma \ref{sigmalemma}. Since $s_n^2/n$ is bounded above and below, the Lindeberg condition with $\sqrt n$ in place of $s_n$ is equivalent to the usual one. We can conclude the Lindeberg Feller condition holds and apply this CLT to \eqref{CF2}, with the justification of the denominator being seen analogously as mentioned above, giving us
		\begin{align}
			&\E_{n,\SC} ( e^{it\s_n^{-1}\langle \lambda^\n,{X^{(n)}- \E_{{n,\b}}X}\rangle} ) \\
			= 	&\frac{\int_{\R} \exp( {-\frac{1}{2n}{\sum_{j=1}^{n} (s^2 \Var_{{n,\b}}(V_j(X_j)) + \frac{t^2\hat{\l}_j^2}{\sigma_n^2} \Var_{{n,\b}}(X_j)+ \frac{2st\hat{\l}_j}{\sigma_n} \Cov_{{n,\b}}(X_j,V_j(X_j)) +\frac{s^2}{2\g})}})ds}{\int_{\R}\exp( {\frac{-1}{2n}{\sum_{j=1}^{n} s^2 (\Var_{{n,\b}}(V_j(X_j))+ \frac{1}{2\g})} })ds} 
			\\&+o(1),\label{CFSC}
		\end{align}
        where we have used Lemma \ref{lemma_b0} to ensure the denominator doesn't tend to zero.
		We now begin to evaluate the $s$ integrals in $\eqref{CFSC}$, the main idea is to view the numerator and denominator as Gaussian densities. In the numerator we first complete the square and rewrite the $s$ terms as
		\begin{align}
			&-\sum_{j=1}^{n} \left(s^2 \Var_{{n,\b}}(V_j(X_j)) + \frac{t^2\hat{\l}_j^2}{\sigma_n^2} \Var_{{n,\b}}(X_j)+ \frac{2st\hat{\l}_j}{\sigma_n} \Cov_{{n,\b}}(X_j,V_j(X_j)) + \frac{s^2}{2\g}  \right)\\
			= & -\left(s\sqrt{\sum_{j=1}^{n}\lb \Var_{{n,\b}}(V_j(X_j))+  \frac{1}{2\g}\rb  }+ \frac{\sum_{j=1}^{n}\frac{t\hat{\l}_j}{\sigma_n}\Cov_{{n,\b}}(V_j(X_j), X_j)}{\sqrt{\sum_{j=1}^{n}(\Var_{{n,\b}}\lb V_j(X_j)) +  \frac{1}{2\g}\rb}}\right)^2 \\ 
			-& \frac{t^2}{\sigma_n^2}\lb\sum_{j=1}^{n}{\hat{\l}_j}^2{\Var_{{n,\b}}(X_j) -\frac{(\sum_{j=1}^{n}\hat{\l}_j\Cov_{{n,\b}}(X_j,V_j(X_j)))^2}{\sum_{j=1}^{n}(\Var_{{n,\b}}(V_j(X_j))+  \frac{1}{2\g}) }}\rb.
		\end{align}
		From here, we can pull the terms not involving $s$ outside of the integral in the numerator, to write $\eqref{CFSC}$ as
		\begin{align}
			&\E_{{n,\SC}}\lb e^{it\s_n^{-1}\langle \lambda,{X- \E_{{n,\b}}X}\rangle} \rb \\ 
			&=  \exp\left( {\frac{-t^2}{2\sigma_n^2}}\left(\sum_{j=1}^{n}{\l_j}^2{\Var_{{n,\b}}(X_j) -\frac{\left(\sum_{j=1}^{n}\l_j\Cov_{{n,\b}}(X_j,V_j(X_j))\right)^2}{\left(\sum_{j=1}^{n}(\Var_{{n,\b}}(V_j(X_j)) +\frac{1}{2\g})\right)}}\right)\right) \\
			&\times\frac{\int_{\R} \exp\left(  {\frac{-1}{2n}\left(s\sqrt{\sum_{j=1}^{n}(\Var_{{n,\b}}(V_j(X_j)) +\frac{1}{2\g} )  }+ \frac{\sum_{j=1}^{n}\frac{t\hat{\l}_j}{\sigma_n}\Cov_{{n,\b}}(V_j(X_j), X_j)}{\sqrt{\sum_{j=1}^{n}(\Var_{{n,\b}}(V_j(X_j)) +\frac{1}{2\g})}}\right)^2}\right)ds} {\left(\int_{\R}\exp\left( {-\frac{1}{2n} s^2 \sum_{j=1}^{n}\left(\Var_{{n,\b}}(V_j(X_j))+\frac{1}{2\g}\right)}\right)ds \right)} + o(1)
			\label{exp}.
		\end{align}
		We now have an exponential term times the product of two $s$ integrals and upon examination, the $s$ integrals in the numerator and denominator of $\eqref{exp}$ can both be seen to be Gaussian integrals with variance given by 
		\begin{align}
			\frac{n}{\sum_{j=1}^{n}\left(\Var_{{n,\b}}(V_j(X_j))+\frac{1}{2\g}\right) }.
		\end{align} 
		Hence, both numerator and denominator can be integrated to give the same normalizing constant which in turn cancel, leaving us with
		\begin{align}
			& \E_{{n,\SC}}\lb e^{it\s_n^{-1}\langle \lambda,{X- \E_{{n,\b}}X}\rangle}\rb \\ 
			&=\exp\left( {\frac{-t^2}{2\sigma_n^2}\left(\sum_{j=1}^{n}{\l_j}^2{\Var_{{n,\b}}(X_j) -\frac{\left(\sum_{j=1}^{n}\l_j\Cov_{{n,\b}}(X_j,V_j(X_j))\right)^2}{\sum_{j=1}^{n}(\Var_{{n,\b}}(V_j(X_j))+\frac{1}{2\g}) }}\right)}\right) +o(1). \label{QQ}
		\end{align}
		From this position we can use the definition of $\sigma_n$, as in Theorem $\ref{T2}$, to see that
		\begin{align}
			\lim_{n \to \infty} \E_{{n,\SC}} e^{it\s_n^{-1}\langle \lambda^\n,{X^{(n)}- \E_{{n,\b}}X}\rangle} = e^{-\frac{t^2}{2}}.
		\end{align}
		We have thus shown that the characteristic function of $P_{n,\SC}$ converges to that of a normal distribution completing the proof.
	\end{proof}

\bigskip
We now give the proof of the hard constraint case, Theorem $\ref{T1}$. This relies upon the proof of the soft constraint case, which we have just developed. 
\begin{proof}[Proof of Theorem \ref{T1}]
	As in Theorem $\ref{T2}$, the proof relies on the key structure of $P_{n}$ developed in Proposition \ref{prop_CF} given by
	\begin{align}
		\E_{n}\lb e^{it\s_n^{-1} \langle \l,X-\E_{n,\b} X \rangle}\rb = \frac{\int_{\R}\E_{{n,\b}}\lb e^{is\lb \frac1{\sqrt n}\sum_{j=1}^n V_j(X_j)-a\sqrt n\rb + it\s_n^{-1}\langle \l,{X-\E_{n,\b} X} \rangle}\rb ds}{\int_{\R}\E_{{n,\b}}\lb e^{is\lb \frac1{\sqrt n}\sum_{j=1}^n V_j(X_j)-a\sqrt n\rb}\rb ds},  \label{CF}
	\end{align}
    which is valid for $n$ sufficiently large. Again we wish to bring the limit $n \to \infty$ inside the $s$ integrals, evaluate this limit using the Lindeberg Feller CLT for a fixed $s,t$ and compute the integral. We will see that taking the limit $n \to \infty$ and evaluating the limiting integrals is analogous to the proof of Theorem $\ref{T2}$. However, the main difficulty of this proof is in bounding the $s$ integrals in $\eqref{CF}$, as unlike the soft constraint case we do not have the Gaussian decay in $s$, see \eqref{Gausdecay}.  We wish to relate the function $\phi$ of Proposition $\ref{lemma_CF_unif}$ to the numerator of $\eqref{CF}$. Here $\hat t_j =t\l_j /\s_n$ and 
	\begin{align}
		||\hat t||_2 = \frac{|t|||\l||_2}{\s_n} \leq \frac{|t|}{\sqrt \e} < \infty,
	\end{align}
	using $\sigma_n^2 \geq \e ||\l||^2$ from the proof  of Lemma \ref{sigmalemma}. Now we can apply Proposition $\ref{lemma_CF_unif}$ for each fixed $t$ once $n$ is sufficiently large, which allows us to write,
	\begin{align}
		&\lvv\E_{{n,\b}}\lb e^{is\lb \frac1{\sqrt n}\sum_{j=1}^n V_j(X_j)-a\sqrt n\rb + it\s_n^{-1}\langle \l,{X-\E_{n,\b} X} \rangle}\rb\rvv= \prod_{j=1}^n \lvv \phi_j\lb\frac{s}{\sqrt{n}}\rb\rvv \\
		& \leq  C'\lb 1+c\min\lb\frac{|s|}{\sqrt n},\frac{s^2}{n}\rb\rb ^{-Cn} 
	\end{align}
	for some $C'>0$, $C>0$ and $c >0$ depending only on $t$. We bound the right hand side by an $n$-independent integrable function, distinguishing two ranges. For $|s| \leq \sqrt n$ the minimum is attained at ${s^2}/{n}$, and here we write
	\begin{align}
		C'\lb 1+c\frac{s^2}{n}\rb^{-Cn} = C' \exp\lb -Cn \log(1+\frac{cs^2}{n})\rb \leq C' \exp\lb -\frac{Cc}{1+c}s^2\rb.
	\end{align}
	For $|s| \geq \sqrt n$ the minimum is attained at ${|s|}/\sqrt n \geq 1$. We assume $s \neq 0$ and split this further into two more cases beginning with $\frac{c^2s^2}{n} \leq |s|$ which gives
	\begin{align}
		&C'\lb 1+\frac{c|s|}{\sqrt n}\rb^{-Cn} \leq C'\lb 1+\frac{c^2s^2}{ n}\rb^{-Cn/2} =  C'\exp\lb-\frac{Cn}{2} \log(1+\frac{c^2s^2}{n})\rb\\
		&\leq  C' \exp\lb-\frac{Cn}{2}\frac{c^2s^2}{(1+|s|)n}\rb = C' \exp\lb-\frac{Cc^2s^2}{2(1+|s|)}\rb.
	\end{align}
	Now the case of $\frac{c^2s^2}{n} \geq |s|$ gives
	\begin{align}
		&C'\lb 1+\frac{c|s|}{\sqrt n}\rb^{-Cn} \leq C'\lb 1+\frac{c^2s^2}{ n}\rb^{-Cn/2} \leq C'\lb \frac{c^2s^2}{ n} \rb ^{-Cn/2} \\
		&= C'\lb \frac{ n}{c^2s^2} \rb ^{Cn/2} \leq C' \lb\frac{1}{|s|} \rb^{Cn/2} = C' {|s|}^{-Cn/2}.
	\end{align}
	Hence the function 
	\begin{align}
		g(s) : = C' \max\lb \exp\lb -\frac{Cc}{1+c}s^2\rb,\exp\lb-\frac{Cc^2s^2}{2(1+|s|)}\rb,  {|s|}^{-\frac{Cn_1}2}\dopp 1_{\{|s|\geq 1\}}\rb,
	\end{align}
	dominates the integrand for all $n \geq n_1$, where $n_1 \geq c^2$ is chosen with $Cn_1 > 2$, and is integrable on $\R$. Lebesgue's dominated convergence theorem permits interchanging the limit $n \to \infty$ with the $s$ integral in the numerator of $\eqref{CF}$, the denominator follows by setting $t=0$. Now, to evaluate
	\begin{align}
		\lim_{n\to\infty} \E_{{n,\b}}\lb e^{is\lb \frac1{\sqrt n}\sum_{j=1}^n V_j(X_j)-a\sqrt n\rb + it\s_n^{-1}\langle \l,{X-\E_{n,\b} X} \rangle}\rb,
	\end{align}
	we define 
	\begin{align}
		Y_j := sV_j(X_j) + \frac{t\hat{\l}_j X_j}{\sigma_n} - \E_{n,\b} \lb sV_j(X_j) + \frac{t\hat{\l}_j X_j}{\sigma_n}\rb,
	\end{align}
	allowing us to re-write this as a sum of independent random variables, giving 
	\begin{align}
		\lim_{n\to\infty} \E_{{n,\b}}\lb e^{is\lb \frac1{\sqrt n}\sum_{j=1}^n V_j(X_j)-a\sqrt n\rb + it\langle \l,\frac{X-\E_{n,\b} X}{\sigma_n} \rangle}\rb = \lim_{n\to\infty} \E_{{n,\b}} \left( e^{\frac{i}{\sqrt{n}}\sum_{j=1}^{n}Y_j}\right).
	\end{align}	
	From this position we can follow the same steps as in the proof of Theorem $\ref{T2}$ to finish the proof, dropping the $1/(2\g)$ term throughout.
\end{proof}

\textbf{Acknowledgments:} The authors are grateful to Fabrice Gamboa for valuable discussions and helpful suggestions.

\printbibliography
\end{document}